\input ./style/arxiv-general.cfg
\documentclass[aos,MSNbibl,nameyear,rotating,dvips]{arximspdf}
\makeatletter
   \@ifpackageloaded{graphicx}{}{\usepackage{graphicx}}
\makeatother

\doi{10.1214/15-AOS1330}
\volume{43}
\issue{5}
\pubyear{2015}
\firstpage{1929}
\lastpage{1958}
\docsubty{FLA}

\makeatletter
\newcommand{\lleft}{\left}
\newcommand{\rright}{\right}
\newtheorem{Th}{Theorem}
\newtheorem{Pro}{Proposition}
\makeatother

\begin{document}
\begin{frontmatter}

\title{Fused kernel-spline smoothing for repeatedly measured
outcomes in a generalized partially linear model with functional
single index\thanksref{T1}}
\runtitle{Fused smoothing for correlated data in single index model}

\begin{aug}
\author[A]{\fnms{Fei}~\snm{Jiang}\thanksref{m1,m2}\corref{}\ead[label=e1]{homebovine@gmail.com}},
\author[B]{\fnms{Yanyuan}~\snm{Ma}\thanksref{m2}\ead[label=e2]{yanyuanma@stat.sc.edu}}
\and
\author[C]{\fnms{Yuanjia}~\snm{Wang}\thanksref{m3}\ead[label=e3]{yuanjiaw@gmail.com}}
\runauthor{F. Jiang, Y. Ma and Y. Wang}
\affiliation{Harvard University\thanksmark{m1}, University of South Carolina\thanksmark{m2} and
Columbia University\thanksmark{m3}}
\address[A]{F. Jiang\\
Harvard University\\
1000 Escalon Avenue\\
Apt. J-1079\\
Sunnyvale, California 94085\\
USA\\
\printead{e1}}
\address[B]{Y. Ma\\
University of South Carolina\\
1523 Greene Street\\
Columbia, South Carolina 29208\\
USA\\
\printead{e2}}
\address[C]{Y. Wang\\
Columbia University\\
722 West 168th St. Rm 205\\
New York, New York 10032\\
USA\\
\printead{e3}}
\end{aug}
\thankstext{T1}{Supported by
NSF Grant DMS-12-06693 and the National Institute of
Neurological Disorders and Stroke Grant NS073671, NS082062.}

%
\received{\smonth{11} \syear{2014}}
%
\revised{\smonth{3} \syear{2015}}

%
\begin{abstract}
We propose a generalized partially linear functional single index
risk score model for repeatedly measured outcomes where the index
itself is a function of time. We fuse the nonparametric kernel method
and regression spline method, and modify the generalized estimating
equation to facilitate estimation and inference.
We use local smoothing kernel to estimate the
unspecified coefficient functions of time, and use
B-splines to estimate the unspecified function of the single index
component. The covariance structure is taken into account via a working
model, which provides valid estimation and inference procedure whether
or not it captures the true covariance. The estimation method is
applicable to both continuous and discrete outcomes.
We derive large sample properties of the estimation procedure and
show a different convergence rate for each component of the model. The
asymptotic properties when the kernel and regression spline methods are
combined in a nested fashion has not been studied prior to this work,
even in the independent data case.
\end{abstract}

\end{frontmatter}

\section{Introduction}\label{secintro}

As a semiparametric regression model, the single index
model is a popular way to accommodate multivariate covariates
while retaining model flexibility. For independent
outcomes, \citet{Carroll1997} introduced a generalized partially linear
single index model, which enriches the family of single index
models by allowing an additional linear component. The goal of this
paper is
to develop a class of generalized partially linear
single index models with functional covariate effect and
explore the estimation and inference for
repeatedly measured dependent outcomes.

In the longitudinal data framework, let $i$
denote the $i$th individual, and $k$ be the $k$th
measurement, where $i=1,\ldots, n$ and $k=1, \ldots, M_i$. Here $M_i$ is the
total number of observations available for the $i$th individual.
Let $D_{ik}$ be the response variable, $\mathbf{Z}_{ik}$ and $\mathbf{X}_{ik}$
be $d_w$- and $d_\beta$-dimensional covariate
vectors.
We assume
the observations from different individuals are independent, while
the responses $D_{i1}, \ldots, D_{iM_i}$, assessed on the same
individual at different time points,
are correlated, but we do not attempt to model such
correlations.
To model the relationship between the conditional mean of the
repeatedly measured outcomes
$D_{ik}$ at time $T_{ik}$
and covariates $\mathbf{Z}_{ik}, \mathbf{X}_{ik}$,
we propose a partially linear functional single index model which models
the mean of $D_{ik}$ given $\mathbf{Z}_{ik}, \mathbf{X}_{ik}$ at time
$T_{ik}$ in the
form of
%
\begin{eqnarray}
\label{eqmodel}
&& E(D_{ik}\vert \mathbf{X}_{ik},
\mathbf{Z}_{ik}, T_{ik})=H\bigl[m\bigl\{\mathbf{w}(T_{ik})^{\mathrm{T}}
\mathbf{Z} _{ik}\bigr\}+\bolds{\beta}^{\mathrm{T}}\mathbf{X}_{ik}
\bigr],
\end{eqnarray}
where $H$ is a known differentiable monotone link function, $\mathbf{w}(t)\in
\mathrm{R}^{d_w}$ at any $t$, $\bolds{\beta}\in
\mathrm{R}^{d_\beta}$. Such a model is useful when the time varying
effect of $\mathbf{Z}_{ik}$ and the
functional combined score effect of $\mathbf{w}(T_{ik})^{\mathrm
{T}}\mathbf{Z}_{ik}$,
adjusted by the covariate vector $\mathbf{X}_{ik}$, are of main
interest.
Note that both
$\mathbf{X}_{ik}$ and $\mathbf{Z}_{ik}$ can contain components that
do not vary with $k$,
such as gender, and the ones
that vary with $k$ such as age. Here, $m(0)$ serves as the
intercept term; thus $\mathbf{X}_{ik}$ does not contain the constant one.
In model~(\ref{eqmodel}), $\mathbf{Z}_{ik}$ includes the covariates
of our main
research interest whose effects are usually time varying and modeled
nonparametrically, and $\mathbf{X}_{ik}$
contains additional covariates of secondary scientific interest and
whose effects are
only modeled via a simple linear form. Here
$m$ is an unspecified smooth
single index function. Further $\mathbf{w}$ is a
$d_w$-dimensional vector of smooth
functions in $L_2$, while $\mathbf{w}(t)$ is $\mathbf{w}$ evaluated
at $t$,
hence a $d_w$-dimensional vector. In addition, $\mathbf{w}(t)$
contributes to
form the argument of
the function $m$, which yields a nested nonparametric
functional form. To ensure identifiability and to reflect the
practical application that motivated this example,
we further
require $\mathbf{w}(t)> \mathbf{0}$ and $\|\mathbf{w}(t)\|_1=
1 \  \forall t$. Here
$\mathbf{w}(t)> \mathbf{0}$ means that every component in $\mathbf
{w}(t)$ is
positive,
and $\|\cdot\|_1$ denotes the vector $l_1$-norm,
that is, the sum of the absolute values of the components in the
vector. The choice of $l_1$ norm incorporates the
practical knowledge from our real data example, described in Section~\ref{secreal}, and is not critical. It can be modified to
other norms, such as the most often used $l_2$ norm or the $\sup$ norm in
our subsequent development. We assume the observed data
follow the model described above. Throughout this paper, we use the
subscript $0$ to denote the
true parameters.
Before we can proceed, we will require the following proposition:

\begin{Pro}\label{proidentify}
Assume $m_0\in\mathcal{M}$, where $\mathcal{M}=
\{m\in
C^1 ([0, 1])$, $m$ is
one-to-one, and $m(0) = c_0\}$.
Here $C^1 ([0, 1])$ is the space of
functions with continuous
derivatives on $[0, 1]$, and $c_0$ is a finite constant.
Assume $\mathbf{w}_{0}(t) \in\mathcal{D}$, where $\mathcal{D} =
\{\mathbf{w}=(w_1, \ldots, w_{d_w})^{\mathrm{T}}\dvtx
\|\mathbf{w}_0(t)\|_1=1, w_j>0
$, and
$w_j\in
C^1 ([0, \tau])\break \forall j=1, \ldots, d_w\}$.
Here
$C^1 ([0, \tau]) $ is the space of
functions with continuous
derivatives on $[0, \tau]$,
and
$\tau$ is a finite constant.
Assume $E(\mathbf{X}_{ik}^{\otimes2})$ and
$E(\mathbf{Z}_{ik}^{\otimes2})$ are both
positive definite, where we define
$\mathbf{a}^{\otimes2} = \mathbf{a}\mathbf{a}^{\mathrm{T}}$ for an arbitrary vector
$\mathbf{a}$.
Then under these assumptions,
the parameter set $(\bolds{\beta}_0,m_0,\mathbf{w}_0)$ in (\ref
{eqmodel}) is identifiable.
\end{Pro}

The proof of Proposition~\ref{proidentify} can be found in
Appendix~\ref{secidentify}.
Model (\ref{eqmodel}) can be viewed as
a longitudinal extension of the generalized partially linear single
index risk
score model introduced in \citet{Carroll1997}; that is,
%
\begin{eqnarray}
\label{eqncarmdl}
&& E(D_{ik}\vert \mathbf{X}_{ik},
\mathbf{Z}_{ik})=H\bigl\{m\bigl(\mathbf{w}^{\mathrm{T}}\mathbf
{Z}_{ik}\bigr)+\bolds{\beta}^{\mathrm{T}}\mathbf{X}_{ik}
\bigr\},
\end{eqnarray}
which is a popular way
to increase flexibility when covariate dimension may be high.
This model was also discussed in Ma et
al. (\citeyear{Ma15}), while the function $m$ was estimated
by using the B-spline technique.
In the
existing literature, many authors explore the generalized partially
linear single index
model under the longitudinal settings. \citet{jiang2011}
consider the single index function in the form of $m(\mathbf{w}^{\mathrm{T}}\mathbf{Z}
_{ik}, t)$, which allows a time
dependent function $m$, but $\mathbf{w}$ is time invariant; hence it
does not
have the nesting structure in model (\ref{eqmodel}) to capture the
time dependent effect of $\mathbf{Z}_{ik}$. Furthermore,
their method does not
consider the within-subject correlation.
\citet{Xu2012} adopted
model (\ref{eqncarmdl}) as a marginal model in the longitudinal data
setting.
Their method takes into account the within-subject correlation, but
similar to the approach of \citet{jiang2011}, it does
not allow $\mathbf{w}$ to vary with time and hence is not sufficient to
describe the
time varying effect of $\mathbf{Z}_{ik}$.
We modify the models of \citet{jiang2011} and \citet{Xu2012} to
accommodate the
time dependent score effect $\mathbf{w}(t)$. In Section~\ref{secreal}, we
show that a time-dependent effect is essential to improve model fit in
some practical situations.
%
%
In addition, we retain the virtue of the models of \citet{jiang2011}
and \citet{Xu2012} by using the semiparametric functional
single index model, which overcomes the curse of dimensionality and
alleviates the risk of
model misspecification [\citep{Peng2011}].

The estimation and inference for model (\ref{eqmodel}) are challenging
due to the nonparametric form of $m, \mathbf{w}$,
and the complications from correlation between
repeatedly measured outcomes. The estimation for
single index models has been discussed extensively in both the
kernel and spline literatures. \citet{Carroll1997} proposed a local
kernel smoothing
technique to estimate the unknown function $m$ and the finite dimensional
parameters $\mathbf{w}, \bolds{\beta}$ in model (\ref{eqncarmdl})
through iterative
procedures. 
Later, \citet{xia2006} applied
a kernel-based minimum
average variance estimation (MAVE)
method for partially linear single index models, which was
first
proposed by \citet{xia2002} for dimension reduction. When $\mathbf{Z}_{ik}$ is
continuous, MAVE results in
consistent estimators for the single index function $m$ without the
root-$n$ assumption on $\mathbf{w}$ as in \citet{Carroll1997}. Nevertheless,
when $\mathbf{Z}_{ik}$ is discrete,
the method may fail to obtain consistent estimators without prior
information about $\bolds{\beta}$
[\citep{xia2002,wang2010}]. Moreover, \citet{wang2009} showed
that MAVE is unreliable for estimating the single index coefficient~$\mathbf{w}$
when $\mathbf{Z}_{ik}$ is unbalanced and sparse, that is, when
$\mathbf{Z}_{ik}$ is
measured at
different time points for each subject, and each subject may have only
a few measurements.

To overcome these limitations,
we apply the B-spline
method to estimate the unknown function $m$, which is stable when
the data set contains discrete or sparse~$\mathbf{Z}_{ik}$. Although the
B-spline method outperforms the kernel method in
estimating~$m$, problems arise if it is also
used for estimating $\mathbf{w}(t)$ in our model setting. If spline
approximations are used for both $m$ and $\mathbf{w}(t)$ with $k$
knots, then
we must simultaneously solve
$(d_w+1)k$ estimating equations to get the
spline coefficients associated with the spline knots,
which may cause numerical
instability and is computationally expensive
when the parameter number increases with the sample size. To
alleviate the computational burden and instability, we estimate~$\mathbf{w}
(t)$ by using
the kernel method. At different time points $t$, the procedure solves
$\mathbf{w}(t)$
independently, and in parallel, it hence does
not suffer from the numerical instability and is computationally efficient.
To handle longitudinal outcomes, we
use the idea from the generalized estimating equation (GEE) to combine
a set of estimating equations built from the
marginal model. It is worth pointing out that the GEE in its
original form is only applicable when the index $\mathbf{w}$ does not change
along time. In conclusion, we combine the kernel and B-spline smoothing
with the
GEE approach, and develop a fused kernel/B-spline procedure
for estimation and inference.

The fusion of
kernel and B-spline poses theoretical challenges which we address in
this work. To the best of our knowledge, this is the first time kernel
and spline methods are jointly implemented in a nested function
setting. We study convergence properties, such as asymptotic bias
and variance, for each component of
the model, show that the parametric component achieves the regular root-$n$
convergence rate and establish the relation of the
nonparametric function convergence rates to the number of B-spline
basis functions and B-spline order, as well as their relation to the
kernel bandwidth. These results provide guidelines for
choosing the number of knots in association with spline order and
bandwidth in order to optimize performance. They also further
facilitate inference, such as constructing confidence intervals and
performing hypothesis testing. Although
theoretical properties of kernel smoothing and spline smoothing are
available separately, the properties, when these two methods are
combined in a nested fashion, have not been studied in the literature,
even for the
independent data case prior to this work. Because the vector function
$\mathbf{w}$ appears inside
the function $m$, the asymptotic analysis of the spline and kernel
methods are not completely separable. This requires a comprehensive
analysis and integration of both methods instead of a mechanical
combination of
two separate techniques.

The rest of the paper is structured as follows.
%
%
In Section~\ref{secest}, we define some notation
and state assumptions in
the model, introduce the fused kernel/B-spline semiparametric
estimating equation, illustrate the
profiling estimation procedure to obtain the estimators and study
the asymptotic properties of the resulting estimators. In Section~\ref{secsim}, we evaluate
the estimation procedure on simulated data sets. In Section~\ref{secreal}, we apply
the model and estimation procedure on the Huntington's disease data
set.
We conclude the paper with some discussion in Section~\ref{seccon}.
We present the technical proofs in the Appendices~\ref{secidentify} and~\ref{secnotationEST}  and an online
supplementary document [\citep{jiang2015}].

\section{Estimating equations and profiling
procedure} \label{secest}

In this section, we construct
estimators for $(\bolds{\beta},m,\mathbf{w})$ in model (\ref
{eqmodel}). We first derive
a set of estimating equations, through applying both B-spline
and kernel methods. We then introduce a profiling procedure to implement
the estimation. Finally, we discuss the asymptotic properties of the estimators.

Many estimation procedures
have been developed for the single index risk score model. In
addition to the methods describe in Section~\ref{secintro}, for the
models with uncorrelated responses, \citet{Cui2012} illustrate an
estimating function method
based on the kernel approach for the generalized single index risk score
model. \citet{Ma2013} discuss a doubly robust and efficient estimation
procedure for the single index risk score model with high-dimensional
covariates. \citet{Ma2014} and \citet{Lu2013} propose B-spline methods
for estimating
the unknown regression link functions in single index risk score models.
However, these methods are not adequate
for the parameter estimation in our model.
As shown in (\ref{eqmodel}), in addition to an unknown link function
$m$, our functional single index model
contains a nonparametric function $\mathbf{w}(t)$ which
is multivariate and appears inside $m$. Therefore, we develop a GEE-type
method for
the parameter estimation in our model, which allows us to take into
account the
within patient correlation. In conjunction
with the kernel smoothing technique and B-spline basis expansion, our
fused method estimates both the
coefficients as a function of time
and the unspecified regression function, and simultaneously handles the
complexities of repeated measurements and curses of dimensionality.

More\vspace*{1pt} specifically, let $\mathbf{B}_r(u) = \{B_{r1}(u),\ldots,
B_{rd_\lambda}(u)\}^{\mathrm{T}}$ be the set of
B-spline basis functions of order $r$,\vspace*{1pt}
and let $\bolds{\lambda}= (\lambda_{1}, \ldots, \lambda_{d_\lambda
})^{\mathrm{T}}
$ be the
coefficients of the B-spline approximation. Denoting $\widetilde{m}(u,
\bolds{\lambda})=\mathbf{B}_r(u)^{\mathrm{T}}\bolds{\lambda}$,
\citet{Boor2001} has shown the existence of
a $\bolds{\lambda}_0 \in
\mathrm{R}^{d_\lambda}$ so that $\widetilde{m}(u,
\bolds{\lambda}_0)=\mathbf{B}_r(u)^{\mathrm{T}}\bolds{\lambda}_0$
converges to $m_0(u)$ uniformly on
$(0, 1)$
when the number of the B-spline inner knots goes to infinity; see \textit{Fact}~1 in Section \textup{S.2} in the supplementary article [\citep{jiang2015}].
A detailed
description of the B-spline functions and the properties of their
derivatives
can be found in \citet{Boor2001}.

The B-spline approximation greatly eases the parameter estimation
procedure. Operationally, for a given sample size $n$, the problem is
reduced from estimating the infinite dimensional $m$ to estimating a
finite dimensional vector
$\bolds{\lambda}$. Since the dimension of $\bolds{\lambda}$ grows
with the sample
size, the estimation consistency can be achieved when the sample size
goes to infinity. Let ${\bolds\theta}=(\bolds{\beta}^{\mathrm
{T}},\bolds{\lambda}^{\mathrm{T}}
)^{\mathrm{T}}\in
\mathbb{R}^{d_\theta}$, and the approximated mean function can be
written as
\begin{eqnarray*}
&& H\bigl[\mathbf{B}_r\bigl\{\mathbf{w}(T_{ik})^{\mathrm{T}}
\mathbf{Z}_{ik}\bigr\} ^{\mathrm{T}}\bolds{\lambda}+\bolds{
\beta}^{\mathrm{T}}\mathbf{X}_{ik}\bigr].
\end{eqnarray*}

We investigate the properties for
estimating $m_0, \mathbf{w}_0, \bolds{\beta}_0$ through
investigating the properties of the estimators for $\bolds{\lambda
}_0$, $\mathbf{w}_0$ and
$\bolds{\beta}_0$.

\subsection{Notation}\label{secnotation}

We define some notation to present the
estimation procedure. To keep the main text
concise, we illustrate the specific forms of notation in
Appendix~\ref{secnotationEST}. Generally, for a
generic vector valued function $\mathbf{a}$ that
depends on some additional parameters, we use
$\widehat{\mathbf{a}}$ to denote the function with the estimated
parameter values
plugged in. For example, this applies to
$\mathbf{S}_w, \mathbf{S}_\beta, \widehat{\mathbf{S}}_{w}$,
$\widehat{\mathbf{S}}_{\beta}$ in the following
text. The specific forms of $\mathbf{S}_w, \mathbf{S}_\beta,
\widehat{\mathbf{S}}_{w},\widehat{\mathbf{S}}_{\beta}$
are\vspace*{1pt} given in  Appendix~\ref{secnotationEST}. 

In our profiling procedure, we
estimate $\bolds{\lambda}_0$ using $\widehat{\bolds{\lambda}}$,
considered as a functional of
$\bolds{\beta}_0, {\mathbf{w}}_0$. Then we estimate $\mathbf{w}_0$ using
$\widehat{\mathbf{w}}$, considered as a function of $\bolds{\beta}_0$
at different time points. Finally, we estimate
$\bolds{\beta}_0$ using $\widehat{\bolds{\beta}}$. We further define
$T_{ik}, k = 1,
\ldots, M_i, i = 1, \ldots, n$ to be the random measurement times
which are independent of $\mathbf{X}_{ik}, \mathbf{Z}_{ik}, D_{ik}$,
$\mathbf{w}$ to be a
function of $t$ for $t\in[0, \tau]$, where $\tau$ is a
finite constant, and $\widehat{\mathbf{w}}(\bolds{\beta}), \widehat{\mathbf{w}}(\bolds{\beta}, t)$, considered as
functions of $\bolds{\beta}$,
to be the estimators for $\mathbf{w}$ and $\mathbf{w}(t)$, respectively.

Let $\mathbf{Q}_\beta(\mathbf{X}_{ik})
= \mathbf{X}_{ik}$,
$\mathbf{Q}_\lambda\{\mathbf{Z}_{ik};\mathbf{w}(t)\}
= \mathbf{B}_r\{\mathbf{w}(t)^{\mathrm{T}}\mathbf{Z}_{ik}\}$ and
$\mathbf{Q}_w\{\mathbf{Z}_{ik};\break  \bolds{\lambda}, \mathbf{w}(t)\}
= \mathbf{Z}_{ik}\mathbf{B}_r'\{\mathbf{w}(t)^{\mathrm{T}}\mathbf
{Z}_{ik}\}^{\mathrm{T}}\bolds{\lambda}$,
be the partial derivatives of\break 
$\mathbf{B}_{r}\{\mathbf{w}(t)^{\mathrm{T}}\mathbf{Z}_{ik}\}
^{\mathrm{T}}\bolds{\lambda}+ \bolds{\beta}^{\mathrm{T}}\mathbf
{X}_{ik}$ with
respect to $\bolds{\beta}$, $\bolds{\lambda}$, $\mathbf{w}(t)$.
In the sequel, we will frequently use $\mathbf{Q}_{\beta ik}$,
$\mathbf{Q}_{\lambda
ik}\{\mathbf{w}(t)\}$, $\mathbf{Q}_{wik}\{ \bolds{\lambda},
\mathbf{w}(t)\}$ as short forms for $\mathbf{Q}_\beta(\mathbf{X}_{ik})$,
$\mathbf{Q}_\lambda\{\mathbf{Z}_{ik};\mathbf{w}(t)\}$ and $\mathbf
{Q}_w\{\mathbf{Z}_{ik}; \bolds{\lambda}, \mathbf{w}
(t)\}$,
respectively.

In general,
to simplify the notation, we use subscripts to indicate the
observations; that is, for a generic function
$a(\cdot)$, we write $a_i(\cdot) \equiv a(\mathbf{O}_i; \cdot)$, where
$\mathbf{O}_i$ denotes the $i$th observed variables. For example, we write
\begin{eqnarray*}
&& H_{ik}\bigl\{\bolds{\beta}, \bolds{\lambda}, \mathbf{w}(t)\bigr\}
\equiv H\bigl[\mathbf{B}_r\bigl\{\mathbf{w}(t)^{\mathrm
{T}}
\mathbf{Z}_{ik}\bigr\}^{\mathrm{T}}\bolds{\lambda}+\bolds{\beta}
^{\mathrm{T}}\mathbf{X}_{ik}\bigr].
\end{eqnarray*}
Further, we indicate the use of the
true function instead of its B-spline approximation by replacing the
argument $\bolds{\lambda}$ with $m$, for example,
\begin{eqnarray*}
&& H_{ik}\bigl\{\bolds{\beta}, m, \mathbf{w}(t)\bigr\}\equiv H\bigl[m
\bigl\{\mathbf{w}(t)^{\mathrm{T}}\mathbf{Z}_{ik}\bigr\} +\bolds{\beta
}^{\mathrm{T}}\mathbf{X}_{ik}\bigr].
\end{eqnarray*}
We also define $\Theta(u) =d
H (u)/d u$ and
\[
\Theta_{ik}\bigl\{\bolds{\beta}, \bolds{\lambda}, \mathbf{w}(t)\bigr\}
= \Theta\bigl[\mathbf{B}_r\bigl\{\mathbf{w}(t)^{\mathrm{T}}
\mathbf{Z}_{ik}\bigr\} ^{\mathrm{T}}\bolds{\lambda}+\bolds{
\beta}^{\mathrm{T}}\mathbf{X}_{ik}\bigr]
\]
and
\[
\Theta_{ik}\bigl\{\bolds{\beta}, m, \mathbf{w}(t)\bigr\} = \Theta
\bigl[m\bigl\{\mathbf{w}(t)^{\mathrm{T}}\mathbf {Z}_{ik}\bigr\} +
\bolds{\beta}^{\mathrm{T}}\mathbf{X}_{ik}\bigr]
\]
throughout the text.

The profiling procedure has three steps. We define the details of the
notation used in
each step and the corresponding population forms in Appendix~\ref{secnotationEST}.

\subsection{Estimation procedure via profiling}\label{secprofile}

In this section, we define the estimation procedures for $m$,
$\mathbf{w}_0$ and $\bolds{\beta}_0$ via estimating equations which
are solved through a
profiling procedure as we describe below. We first
estimate the function $m$ through B-splines,
by treating $\mathbf{w}$ and $\bolds{\beta}$ as parameters that are
held fixed. This
yields a set of estimating equations for the spline coefficients, as
functions of $\mathbf{w}$ and $\bolds{\beta}$.
We then estimate the partially linear nonparametric component $\mathbf{w}(t)$
of the
cognitive score profiles through local kernel smoothing,
while treating $\bolds{\beta}$ as fixed parameters. This further
allows us to
obtain a second set of estimating equations at each time point that the
function $\mathbf{w}(t)$ needs to be estimated, as a function of~$\bolds{\beta}$.
Finally, we estimate the parametric component coefficients $\bolds
{\beta}$
through solving its own corresponding estimating equation set.
The profiling procedure achieves a certain separation by allowing us
to treat only one of the three components in
each of the three nested steps; hence it
eases the computational complexities. Because the
B-spline estimator $\widehat{\bolds{\lambda}}$, kernel estimator
$\widehat{\mathbf{w}}(t)$ and linear parametric estimator
$\widehat{\bolds{\beta}}$ have different convergence rates, such separation
also facilitates analysis of the asymptotic properties,
compared with a simultaneous estimation procedure.

\begin{longlist}[{}]
\item[\textit{Step} 1.]
We obtain $\widehat{\bolds{\lambda}}(\bolds{\beta}_0, {\mathbf
{w}}_0)$ by solving
\begin{eqnarray}
{\label{eqnstep1}}
&&\sum_{i=1}^{n}\widetilde{
\mathbf{Q}}_{\lambda i}\bigl\{\mathbf {w}_0(\mathbf{T}_i)
\bigr\} ^{\mathrm{T}}\bolds{\Theta}_i \bigl\{\bolds{\beta}
_0, \bolds{\lambda}, \mathbf{w}_0(\mathbf{T}_i)
\bigr\} \bolds{\Omega}^{-1}_i \bigl[\mathbf{D}_i
- \mathbf{H}_i \bigl\{\bolds{\beta }_0, \bolds{\lambda},
\mathbf{w}_0(\mathbf{T}_i)\bigr\}\bigr] = {\mathbf0}
\nonumber
\end{eqnarray}
with respect to $\bolds{\lambda}$,
where
$\bolds{\Omega}_i$ is a working covariance matrix, and
$
\bolds{\Theta}_i = \operatorname{diag}\{
\Theta_{ik}\}$, $k = 1, \ldots, M_i$
is a $M_i \times M_i$ diagonal matrix. From the first step, we obtain
the B-spline coefficients to estimate the
function $m$.

\item[\textit{Step} 2.]
We obtain $\widehat{\mathbf{w}}(\bolds{\beta})$ in this step.
Let $\mathbf{K}_h(\mathbf{T}_i - t_0)$ be a $d_w M_i \times d_w M_i$
diagonal matrix whose $k$th diagonal block is
$
\operatorname{diag}\{K_h(T_{ik} - t_0)\}
$
where $K_h(s) = h^{-1} K(s/h)$ is a Kernel function with
bandwidth $h$.

To obtain $\widehat{\mathbf{w}}(\bolds{\beta}_0, t_0)$, we
solve the estimating equation
%
\begin{eqnarray}
&&\sum_{i=1}^{n}\widehat{
\mathbf{A}}_{wi}\bigl\{\bolds{\beta}_0, \widehat{\bolds{
\lambda}}(\bolds{\beta}_0, {\mathbf{w}}),\mathbf{w}(t_0)
\bigr\}\widehat{\mathbf{V}}_{wi}\bigl\{\bolds {\beta}_0,
\widehat{\bolds{\lambda}}(\bolds{\beta}_0, {\mathbf{w}}),
\mathbf{w}(t_0)\bigr\}^{-1}
\nonumber
\\[-8pt]
\label{eqnstep2}
\\[-8pt]
\nonumber
&&\hspace*{2pt}\quad{}\times \mathbf{K}_h(\mathbf{T}_i - t_0)
\widehat{\mathbf{S}}_{wi}\bigl\{\bolds{\beta}_0, \widehat{
\bolds {\lambda}}(\bolds{\beta}_0, {\mathbf{w}}),\mathbf{w}(t_0)
\bigr\}
\end{eqnarray}
with respect to ${\mathbf{w}}$.
Recall that $\|\mathbf{w}(t_0)\|_1 =
1$.
In the implementation,
we parameterize $w_{d_w}= 1 -\sum_{j=1}^{d_w-1} w_j$, and
derive the score functions for the vector $(w_1, \dots, w_{d_w-1})$.
We then solve the estimating equation system which contains the $d_w-1$
equations
constructed from the score
functions and the equation $\sum_{j=1}^{d_w} w_j - 1 = 0$. The roots\vspace*{1pt}
of the estimating equation system automatically satisfy the $l_1$ constraint.
In all our experiments, the resulting $\widehat{w}_j(t)$ are nonnegative
automatically, and hence we did not particularly enforce the nonnegativity
as a constraint. If it is needed, one can further enforce the
nonnegativity and perform a constrained optimization.

\item[\textit{Step} 3.]
We obtain $\widehat{\bolds{\beta}}$ by solving
%
\begin{eqnarray}
&&\sum_{i=1}^{n}\widehat{\mathbf{A}}_{\beta i}\bigl[\bolds{\beta}, \widehat{\bolds{\lambda}}\bigl\{
\bolds{\beta}, \widehat{\mathbf {w}}(\bolds{\beta}, \mathbf{T}_{i})
\bigr\}, {\widehat{\mathbf{w}}}(\bolds{\beta})\bigr]\widehat{\mathbf{V}}_{\beta
i}
\bigl[\bolds{\beta}, \widehat{\bolds{\lambda}}\bigl\{\bolds{\beta}, {\widehat{
\mathbf{w}}}(\bolds{\beta})\bigr\}, \widehat{\mathbf {w}}(\bolds{\beta},
\mathbf{T}_{i})\bigr]^{-1}
\nonumber
\\[-8pt]
\label{eqnstep3}
\\[-8pt]
\nonumber
&&\hspace*{3pt}\quad{}\times\widehat{\mathbf{S}}_{\beta i}\bigl[\bolds{\beta}, \widehat{\bolds{
\lambda}}\bigl\{\bolds{\beta}, {\widehat{\mathbf{w}}}(\bolds {\beta} )\bigr\},
\widehat{\mathbf{w}}(\bolds{\beta}, \mathbf{T}_{i})\bigr]=\mathbf{0}.
\nonumber
\end{eqnarray}

In above steps, we approximate $\partial\widehat{\mathbf{w}}(\bolds
{\beta}, \mathbf{T}
_{i})/\partial
\bolds{\beta}^{\mathrm{T}}$, $\partial\widehat{\bolds{\lambda
}}(\bolds{\beta}, \mathbf{w})/\partial\bolds{\beta}^{\mathrm{T}}$
and\break $\partial\widehat{\bolds{\lambda}}(\bolds{\beta}_0,
{\mathbf{w}}_0)/\partial{\mathbf{w}}$ by
the leading terms in their expansions. Their explicit forms are shown in
(S.27) in the proofs of Lemma~6,
(S.37) in the proofs of Lemma~11 and the notation in  \textit{step}~2 in
the Appendix~\ref{secnotationEST}, respectively.
\end{longlist}

\subsection{Asymptotic properties of the estimators}\label{secthm}

The profiling estimator described in Section~\ref{secprofile}
is quite complex, caused by the functional nature of $\mathbf{w}(t)$,
the unspecified forms of
both $\mathbf{w}$ and $m$ and their nested appearance in the model,
the correlation among different observations associated with the same
individual and the different numbers of observations for each individual.
In addition, the fused kernel/B-spline method requires careful
joint consideration of both smoothing techniques. As a consequence, the analysis
to obtain the asymptotic properties of the estimator described in
Section~\ref{secprofile} is very challenging and involved.
We first list the regularity conditions under which we
perform our theoretical analysis:
\begin{longlist}[(A6)]
\item[(A1)] The kernel function $K(\cdot)$ is nonnegative,
has compact support
and satisfies $\int K(s) \,ds=1$, $\int K(s)s \,ds=0$ and $\int K(s)
s^2 \,ds < \infty$, and\break  $\int K^2(s)
s \,ds < \infty$.

\item[(A2)] The bandwidth $h$ in the kernel smoothing satisfies
$nh^2\rightarrow\infty$ and $nh^4
\rightarrow0$ when $n\to\infty$.

\item[(A3)] The density function of $\mathbf{w}(t)^{\mathrm
{T}}\mathbf{Z}$ for each $t \in
[0, \tau]$ is bounded away
from~0 on $S_w(t)$ and satisfies the Lipschitz condition of order 1
on $S_w(t)$, where $\mathbf{w}$ is in a neighborhood of $\mathbf
{w}_0$, and $S_w(t) =
\{\mathbf{w}(t)^{\mathrm{T}}\mathbf{Z}, \mathbf{Z}\in S \}$ and $S$
is a
compact support of $\mathbf{Z}$ and $\tau<\infty$ is a finite constant.
Without loss of generality, we assume $S_w(t)
= [0, 1]$.

\item[(A4)] Assume $m_0\in
\{m\in
C^q ([0, 1])$, $m$ is
one-to-one and $m(0) = c_0\} $. Here $C^q ([0, 1]) $ is the space of
functions with first $q$ continuous
derivatives on $[0, 1]$.
The spline order is $r \geq
q$. The cluster
size $M_i$ is a fixed finite number that does
not diverge with the sample size, that is, $M_i < \infty$ for all $i$.
\item[(A5)] Let $h_p$ be the distance between the $(p+1)$th and $p$th
interior knots of the order $r$ B-spline
functions. And $h_b = \max_{r\leq p\leq N+r} h_{p}$. There
exists $0 <c_{h_b} < \infty$, such that
\begin{eqnarray*}
&& \max_{r \leq p \leq N+r} h_{ p+1}= o\bigl(N^{-1}\bigr)
\quad\mbox{and}\quad h_b / \min_{r \leq p \leq N+r} h_{p} <
c_{h_b},
\end{eqnarray*}
where $N$ is the number of knots which satisfies $N \rightarrow
\infty$ as $n\rightarrow\infty$, and $N^{-1} n(\log n )^{-1}
\rightarrow\infty$ and $N n^{-1/(2 q + 1)} \rightarrow\infty$,
further assuming $q > 3$ and $N^{-3} n \to\infty$.
\item[(A6)]
The matrices\vspace*{1pt}
$E(\mathbf{X}_{ik}^{\otimes2})$,
$E([\mathbf{X}_{ik} - E\{\mathbf{X}_{ik}|\mathbf{w}(t)^{\mathrm
{T}}\mathbf{Z}_{ik}\}]^{\otimes2})$,
$E([\mathbf{Z}_{ik} -  E\{\mathbf{Z}_{ik} |\mathbf{w}(t)^{\mathrm
{T}}\mathbf{Z}_{ik}\}
m'_0\{\mathbf{w}(t)^{\mathrm{T}}\mathbf{Z}_{ik}\}]^{\otimes2})$ and
$E([\mathbf{X}_{ik}\mathbf{Z}
_{ik}^{\mathrm{T}}-
E\{\mathbf{X}_{ik}\mathbf{Z}_{ik}^{\mathrm{T}}|\mathbf
{w}(t)^{\mathrm{T}}\mathbf{Z}_{ik}\}\times\break m'_0\{\mathbf{w}(t)^{\mathrm
{T}}\mathbf{Z}
_{ik}\}]^{\otimes2})$
are finite and positive definite for any $t\in[0,
\tau]$.
\end{longlist}

The requirements $nh^4\to0$ in (A2) and $N n^{-1/(2 q + 1)}
\rightarrow \infty$ in (A5) are undersmoothing requirements on the kernel
approximation and on the spline
approximation, respectively. They are required to ensure that
the biases, $E(\widehat{\mathbf{w}})-\mathbf{w}$ and $E(\mathbf
{B}_r^{\mathrm{T}}\widehat{\bolds{\lambda}})-m_0$
are ignorable compared to other terms left in the final analysis.
These kinds of undersmoothing
conditions are commonly required in semiparametric models.

Theorems \ref{th1}--\ref{th3} describe the asymptotic properties
for the
estimators of $\mathbf{w}_0(t)$, $\bolds{\beta}_0$ and $m_0$, respectively.
%
\begin{Th}\label{th1}
Assume\vspace*{1pt} conditions \textup{(A1)}--\textup{(A6)} and the identifiability conditions stated
in Proposition~\ref{proidentify} hold. Let $\widehat{\mathbf{A}}_{wi}$, $\widehat{\mathbf{V}}_{wi}$, $\widehat{\mathbf{S}
}_{wi}$ and their
population forms ${\mathbf{A}}_{wi}$, ${\mathbf{V}}_{wi}$ and
${\mathbf{S}}_{wi}$ be as
defined in the notation of \textup{step 2} in Appendix~\ref{secnotationEST}. Let $ {\widehat{\mathbf{w}}}(\bolds{\beta}_0,
t_0)$ solve (\ref{eqnstep2}) and
$f_T$ be the probability density function of $T_{ik}$ with support
$[0, \tau]$. Define\vspace*{-3pt}
\begin{eqnarray*}
\Sigma_w & =& (nh)^{-1} \bigl\{\mathbf{B}(t_0)
f_T(t_0)\bigr\}^{-1}
\\[-2pt]
&&{}\times E \biggl(f_T(t_0)
\bigl[{\mathbf{A}}_{wi}\bigl\{\bolds{\beta}_0,m_0,
\mathbf {w}_0(t_0)\bigr\}{\mathbf{V}}_{wi}\bigl\{\bolds{\beta}_0,
m_0,\mathbf{w}_0(t_0)\bigr\} ^{-1}
\bigr]\\[-2pt]
&&\hspace*{27pt}{}\times\int\mathbf{K}(\mathbf{s}) \mathbf{V} _{w i}^*\bigl\{\bolds{
\beta}_0, m_0,\mathbf{w}_0(t_0)
\bigr\} \mathbf{K}(\mathbf{s})\, d\mathbf{s}
\\[-2pt]
&&\hspace*{38pt}{}\times\bigl[{\mathbf{A}}_{wi}\bigl\{\bolds{\beta}_0,
m_0,\mathbf{w}_0(t_0)\bigr\}{
\mathbf{V}}_{wi}\bigl\{\bolds{\beta}_0, m_0,
\mathbf{w}_0(t_0)\bigr\}^{-1} \bigr]
^{\mathrm{T}} \biggr)
\\[-2pt]
&&{}\times\bigl\{\mathbf{B}(t_0)f_T(t_0)\bigr
\}^{-1}.
\end{eqnarray*}
Then\vspace*{-6pt}
\begin{eqnarray*}
&& \Sigma_w^{-1/2}\bigl\{\widehat{\mathbf{w}}(\bolds{
\beta}_0, t_0) - \mathbf{w}_0(t_0)
\bigr\} \stackrel{d} {\rightarrow} N(0, \mathbf{I}),
\end{eqnarray*}
where $\mathbf{B}$ are defined in the notation of \textup{step}~3 in Appendix~\ref{secnotationEST}.
\end{Th}

Theorem~\ref{th1} establishes the large sample properties
of the estimation of the multivariate weight function $\mathbf{w}_0(t)$.
It shows that our method achieves the usual
nonparametric convengence rate of root-$nh$ under the conditions
given.

\begin{Th}\label{th2}
Assume\vspace*{1pt} conditions \textup{(A1)}--\textup{(A6)} and the identifiability conditions stated
in Proposition~\ref{proidentify} hold. Let $\widehat{\mathbf{S}}_{\beta i k}$,
$\widehat{\mathbf{A}}_{\beta i k}$, $\widehat{\mathbf{V}}_{\beta i
k l}$ and their population
forms ${\mathbf{S}}_{\beta i k}$,
${\mathbf{A}}_{\beta i k}$, ${\mathbf{V}}_{\beta i k l}$ be as
defined in the
notation of \textup{step} 3 in Appendix~\ref{secnotationEST}, and
${\widehat{\mathbf{w}}}(\bolds{\beta})$,
$ {{\mathbf{w}}}(\bolds{\beta}) $ be as defined in Section~\ref{secnotation}.
Let $\widehat{\bolds{\beta}}$ solve (\ref{eqnstep3}),\vspace*{-2pt}
then
%
\begin{eqnarray}
&&\sqrt{n} (\widehat{\bolds{\beta}}- \bolds{
\beta}_0)
\nonumber
\\[-2pt]
&&\qquad =\mathbf{F}(m_0) ^{-1} \Biggl(\frac{1}{\sqrt{n}}\sum
_{i=1}^{n}{\mathbf{A}}_{\beta
i}\bigl\{
\bolds{\beta}_0, m_0, \mathbf{w}_0(
\mathbf{T}_i)\bigr\}{\mathbf {V}}_{\beta
i}\bigl\{\bolds{
\beta}_0, m_0, \mathbf{w}_0(
\mathbf{T}_i)\bigr\} ^{-1}
\nonumber
\\[-2pt]
\nonumber
&&\hspace*{80pt}\qquad\quad{}\times {\mathbf{S}}_{\beta i}\bigl\{\bolds{\beta}_0,
m_0, \mathbf{w}_0(\mathbf{T}_i)\bigr\}
\\[-2pt]
\nonumber
&&\hspace*{45pt}\qquad\quad{} -
\frac{1}{\sqrt{n}}\sum_{j=1}^{n}E \bigl(
\mathbf{A}_{\beta
i}\bigl\{\bolds{\beta}_0, m_0,
\mathbf{w}_0(\mathbf{T}_j)\bigr\}\\
&&\hspace*{93pt}\nonumber\qquad\qquad{}\times
\mathbf{V}_{\beta
i}\bigl\{\bolds{\beta}_0,
m_0, \mathbf{w}_0(\mathbf{T}_i)\bigr
\}^{-1}\mathcal{K}(\mathbf{T}_j) | \mathbf{O}_j
\bigr)\mathbf{B}(\mathbf{T}_{j})^{-1} \\[-2pt]
\nonumber
&&\hspace*{82pt}\qquad\qquad{}\times\bigl[ {
\mathbf{A}}_{wj}\bigl\{\bolds{\beta}_0, m_0,
\mathbf{w}_0(\mathbf {T}_{j})\bigr\}
\\[-2pt]
\nonumber
&&\hspace*{107pt}\qquad\quad{}\times{\mathbf{V}}_{wj}\bigl\{\bolds{\beta}_0,
m_0,\mathbf{w}_0(\mathbf{T}_{j})\bigr
\}^{-1}\\[-2pt]
\label{eqnbetaexp}
&&\hspace*{198pt}\qquad\quad{}\times {\mathbf{S}}_{wj}\bigl\{\bolds{\beta}_0,
m_0,\mathbf{w}_0(\mathbf {T}_{j})\bigr\}
\bigr] \Biggr)
\\[-2pt]
&&\qquad\quad{}-\mathbf{G}(m_0) \mathbf{V}^{-1}\frac{1}{\sqrt{n}}\sum_{i=1}^{n}
\widetilde{\mathbf {Q}}_{\lambda i}\bigl\{ \mathbf{w}_0(\mathbf{T}
_i)\bigr\} ^{\mathrm{T}}\bolds{\Theta}_i\bigl\{\bolds{
\beta}_0, m_0, \mathbf{w}_0(
\mathbf{T}_i)\bigr\}
\nonumber
\\[-2pt]
&&\hspace*{94pt}\qquad\quad{}\times\bolds{\Omega}^{-1}_i
\bigl[\mathbf{D}_i- \mathbf{H}_i \bigl\{\bolds{\beta}_0,
m_0, \mathbf{w}_0(\mathbf{T}_i)\bigr\}\bigr]
 \bigl\{1 + o_p(1)\bigr\},
\nonumber
\end{eqnarray}
where $\mathcal{K}(\mathbf{T}_i )= \operatorname{diag}\{\bolds
{\kappa}(T_{ik}), k = 1, \ldots, M_i\}$
a $d_\beta M_i \times d_w M_i$ matrix
and $\bolds{\kappa}(T_{ik})$ is
\begin{eqnarray*}
&& \biggl\{\mathbf{Q}_{\beta i k} - \bolds{\delta}\bigl\{\mathbf
{w}_0(T_{ik})^{\mathrm{T}}\mathbf{Z}_{ik}\bigr\}
\\
&&\hspace*{6pt}{}- \biggl({\mathbf{B}}(T_{ik})^{-1} E \biggl[{
\mathbf{A}}_{wj}\bigl\{\bolds {\beta}_0, m_0,
\mathbf{w}_0(T_{ik})\bigr\}{\mathbf{V}}_{wj}\bigl\{\bolds{\beta}_0,
m_0,\mathbf{w}_0(T_{ik})\bigr\} ^{-1}
\\
&&\hspace*{110pt}\qquad\qquad{}\times\frac{\partial{\mathbf{S}
}_{wj}\{\bolds{\beta}_0,
m_0,\mathbf{w}_0(T_{ik})\} }{\partial\bolds{\beta}^{\mathrm
{T}}}\Big\vert \mathbf{O}_i \biggr] \biggr)^{\mathrm{T}}
\mathbf{Z}_{ik}
\\
&&\hspace*{143pt}{}\times m'_0\bigl\{\mathbf{w}_0(T_{ik})^{\mathrm{T}}
\mathbf{Z}_{ik}\bigr\}+ \bolds{\gamma}\bigl\{\mathbf{w}_0(T_{ik})^{\mathrm{T}}
\mathbf{Z}_{ik}\bigr\} \biggr\}
\\
&&\qquad{}\times\mathbf{Q}_{wik}\bigl\{
m_0,\mathbf{w}_0(T_{ik})\bigr\}^{\mathrm{T}}\Theta_{ik} \bigl\{\bolds{\beta}_0,
m_0,\mathbf{w}_0(T_{ik})\bigr\}.
\\
&& \mathbf{F}(m_0) = -E \biggl\{{\mathbf{A}}_{\beta i}\bigl\{
\bolds{\beta }_0, m_0, \mathbf{w}_0(
\mathbf{T}_{i})\bigr\} {\mathbf{V}}_{\beta
i}\bigl\{\bolds{
\beta}_0, m_0, \mathbf{w}_0(
\mathbf{T}_i)\bigr\}^{-1}
\\
&&\qquad\hspace*{128pt}{}\times\frac{\partial{\mathbf{S}}_{\beta
i}\{\bolds{\beta}_0, m_0, \mathbf{w}_0(\mathbf{T}_i)\}}{\partial
\bolds{\beta}^{\mathrm{T}}} \biggr\}
\end{eqnarray*}
and
\begin{eqnarray*}
\mathbf{G}(m_0) &=& E\bigl[ {\mathbf{A}}_{\beta
i}\bigl\{
\bolds{\beta}_0, m_0, \mathbf{w}_0(
\mathbf{T}_i)\bigr\}{\mathbf{V}}_{\beta
i}\bigl\{\bolds{
\beta}_0, \bolds{\lambda}_0 , \mathbf{w}_0(
\mathbf{T}_i)\bigr\}^{-1} \mathbf{C}_i
\\
&&\qquad\hspace*{35pt}{}\times\bolds{\Theta}_i^*\bigl\{\bolds{\beta}_0,
m_0, {\mathbf{w}}_0(\mathbf{T}_i)\bigr\}
\mathbf{Q}^*_{\lambda i}\bigl\{\mathbf{w}_0(\mathbf{T}_i)
\bigr\}\bigr].
\end{eqnarray*}
Here $\mathbf{C}_i$ is a $d_\beta M_i \times d_\beta M_i$ with the $k$th
block having the form
\begin{eqnarray*}
&& \biggl\{\mathbf{Q}_{\beta i k} - \bolds{\delta}\bigl\{\mathbf
{w}_0(T_{ik})^{\mathrm{T}}\mathbf{Z}_{ik}\bigr\}
\\
&&\hspace*{6pt}{}- \biggl({\mathbf{B}}(T_{ik})^{-1} E \biggl[{
\mathbf{A}}_{wj}\bigl\{\bolds {\beta}_0, m_0,
\mathbf{w}_0(T_{ik})\bigr\}{\mathbf{V}}_{wj}\bigl\{\bolds{\beta}_0,
m_0,\mathbf{w}_0(T_{ik})\bigr\} ^{-1}
\\
&&\hspace*{155pt}{}\times\frac{\partial{\mathbf{S}}_{wj}\{\bolds{\beta}_0,
m_0,\mathbf{w}_0(T_{ik})\} }{\partial\bolds{\beta}^{\mathrm
{T}}}\Big\vert \mathbf{O}_i \biggr]
\biggr)^{\mathrm{T}}
\\
&&\hspace*{70pt}\qquad\qquad{}\times\mathbf{Z}_{ik} m_0'\bigl
\{\mathbf {w}_0(T_{ik})^{\mathrm{T}}\mathbf{Z}_{ik}
\bigr\}+ \bolds{\gamma}\bigl\{\mathbf{w}_0(T_{ik})^{\mathrm{T}}
\mathbf {Z}_{ik}\bigr\} \biggr\}.
\end{eqnarray*}
Here $\bolds{\Theta}_i^*\{\bolds{\beta}_0, m_0,
{\mathbf{w}}_0(\mathbf{T}_i)\}$ is a $d_\beta M_i \times d_\beta M_i$
matrix with the
$k$th block being a $d_\beta\times d_\beta$ diagonal matrix with the
element $\Theta_{ik}\{\bolds{\beta}_0, m_0,
{\mathbf{w}}_0(T_{ik})\}$. And
$\mathbf{Q}^*_{\lambda i}\{\mathbf{w}_0(\mathbf{T}_i)\}$ is a
$d_\beta M_i \times d_\lambda$
matrix with $k$th row block being a $d_\beta\times d_\lambda$
matrix, where $d_\beta$ replicates of the row vector
$\mathbf{Q}_{\lambda i k}\{\mathbf{w}_0(T_{ik})\}^{\mathrm{T}}$.
$\mathbf{B}$, $\bolds{\delta}$, $\bolds{\gamma}$
are functions
defined in the notation of \textup{step}~3 in Appendix~\ref{secnotationEST}.

Consequently,\vspace*{-3pt} we have
\begin{eqnarray*}
&& \sqrt{n} (\widehat{\bolds{\beta}}- \bolds{\beta}_0)\stackrel {d} {
\rightarrow}N(0, \Sigma),
\end{eqnarray*}
where\vspace*{-2pt}
\begin{eqnarray*}
\Sigma&=&\mathbf{F}(m_0) ^{-1} E \bigl[ \bigl(\bigl[{
\mathbf{A}}_{\beta
i}\bigl\{\bolds{\beta}_0, m_0,
\mathbf{w}_0(\mathbf{T}_i)\bigr\}{\mathbf
{V}}_{\beta
i}\bigl\{\bolds{\beta}_0, m_0,
\mathbf{w}_0(\mathbf{T}_i)\bigr\}^{-1}
\\
&&\hspace*{153pt}{}\times {\mathbf{S}}_{\beta i}\bigl\{\bolds{\beta}_0,
m_0, \mathbf{w}_0(\mathbf{T}_i)\bigr\}
\bigr]^{\otimes2} \bigr)\\
&&\hspace*{29pt}\qquad{}+ \bigl\{ \bigl(E \bigl(\mathbf {A}_{\beta
i}
\bigl\{\bolds{\beta}_0, m_0, \mathbf{w}_0(
\mathbf{T}_j)\bigr\}
\mathbf{V}_{\beta
i}\bigl\{\bolds{\beta}_0,
m_0, \mathbf{w}_0(\mathbf{T}_i)\bigr
\}^{-1}\mathcal{K}(\mathbf{T}_j) | \mathbf{O}_j
\bigr)\\
&&\hspace*{29pt}\hspace*{16pt}\qquad{}\times\mathbf{B}(\mathbf{T}_{j})^{-1} \bigl[ {
\mathbf{A}}_{wj}\bigl\{\bolds{\beta}_0, m_0,
\mathbf{w}_0(\mathbf {T}_{j})\bigr\}
\\
&&\hspace*{120pt}{}\times{\mathbf{V}}_{wj}\bigl\{\bolds{\beta}_0,
m_0,\mathbf{w}_0(\mathbf{T}_{j})\bigr
\}^{-1} \\
&&\hspace*{205pt}{}\times{\mathbf{S}}_{wj}\bigl\{\bolds{\beta}_0,
m_0,\mathbf{w}_0(\mathbf {T}_{j})\bigr\}
\bigr] \bigr)^{\otimes2} \bigr\}
\\
&&\hspace*{50pt}{}+ \bigl\{ \bigl( \mathbf{G}(m_0) \mathbf{V}^{-1}
\widetilde{\mathbf {Q}}_{\lambda i}\bigl\{ \mathbf{w}_0(
\mathbf{T}_i)\bigr\} ^{\mathrm{T}}\bolds{\Theta}_i\bigl
\{\bolds{\beta}_0, m_0, \mathbf{w}_0(
\mathbf{T}_i)\bigr\}
\\
&&\hspace*{160pt}{}\times\bolds{\Omega}^{-1}_i
\bigl[\mathbf{D}_i- \mathbf{H}_i \bigl\{\bolds{\beta}_0,
m_0, \mathbf{w}_0(\mathbf{T}_i)\bigr\}\bigr]
\bigr)^{\otimes2} \bigr\} \bigr]
\\
&&{}\times\mathbf{F}(m_0)^{-1}.
\end{eqnarray*}
\end{Th}
Theorem~\ref{th2} establishes the usual parametric convergence rate
for $\widehat{\bolds{\beta}}$, even though the estimation relies on
multiple nonparametric estimates as well.
The form of~(\ref{eqnbetaexp}) in Theorem~\ref{th2} indicates that
the variance of estimating
$\bolds{\beta}_0$ is inflated by the estimation $\widehat{\mathbf
{w}}$, as given in
\begin{eqnarray*}
&&\frac{1}{\sqrt{n}}\sum_{j=1}^{n}E \bigl(
\mathbf{A}_{\beta
i}\bigl\{\bolds{\beta}_0, m_0,
\mathbf{w}_0(\mathbf{T}_j)\bigr\}\mathbf{V}_{\beta
i}
\bigl\{\bolds{\beta}_0, m_0, \mathbf{w}_0(
\mathbf{T}_i)\bigr\}^{-1}\mathcal{K}(\mathbf{T}_j)
| \mathbf{O}_j \bigr)
\\
&&\hspace*{12pt}\qquad{}\times\mathbf{B}(\mathbf{T}_{j})^{-1} \bigl[ {
\mathbf{A}}_{wj}\bigl\{\bolds{\beta}_0, m_0,
\mathbf{w}_0(\mathbf{T}_{j})\bigr\}{\mathbf{V}}_{wj}
\bigl\{\bolds{\beta}_0, m_0,\mathbf{w}_0(
\mathbf{T}_{j})\bigr\}^{-1}
\\
&&\hspace*{112pt}\hspace*{50pt}\qquad{}\times {\mathbf{S}}_{wj}\bigl\{\bolds{\beta}_0,
m_0,\mathbf{w}_0(\mathbf{T}_{j})\bigr\} \bigr]
\end{eqnarray*}
and is also inflated by the estimation
$\widehat{\lambda}$, as given in
\begin{eqnarray*}
&&\mathbf{G}(m_0) \mathbf{V}^{-1} \frac{1}{\sqrt{n}}\sum
_{i=1}^{n}\widetilde{\mathbf{Q}}_{\lambda i}
\bigl\{ \mathbf{w}_0(\mathbf{T}_i)\bigr\} ^{\mathrm{T}}
\bolds{\Theta}_i\bigl\{\bolds {\beta}_0, m_0,
\mathbf{w}_0(\mathbf{T}_i)\bigr\}
\\
&&\hspace*{39pt}\qquad\qquad{}\times \bolds{
\Omega}^{-1}_i \bigl[\mathbf{D}_i- \mathbf{H}_i \bigl\{\bolds{\beta}_0,
m_0, \mathbf{w}_0(\mathbf{T}_i)\bigr\}
\bigr].
\end{eqnarray*}
See Lemmas~9, 11 and the proofs of Theorem~\ref{th2} in the supplementary article [\citep{jiang2015}] for a more
detailed discussion.

The asymptotic normality of $\widehat{\bolds{\beta}}$ established in
Theorem~\ref{th2}
further facilitates inference on $\bolds{\beta}$ such as constructing
confidence
intervals or performing hypothesis testing. In implementing these
inference procedures, we replace the variance--covariance matrix
$\Sigma$ with its estimate, where we use empirical sample mean over
the observed
samples to replace the
expectations\vadjust{\goodbreak} in Theorem~\ref{th2}, and plug in the estimates of
the corresponding parameter and function values. This is the procedure
adopted in all our numerical\vspace*{-2pt} implementation.

\begin{Th}{\label{th3}}
Assume conditions \textup{(A1)}--\textup{(A6)} and the identifiability conditions stated in
Proposition~\ref{proidentify} hold. Let $\widehat{m}\{u,
\widehat{\bolds{\lambda}}(\bolds{\beta}, {\mathbf{w}})\} =
\mathbf{B}_r(u) ^{\mathrm{T}}\widehat{\bolds{\lambda}}(\bolds
{\beta},
{\mathbf{w}})$, $\widetilde{m}\{u,
\bolds{\lambda}_0\} = \mathbf{B}_r(u) ^{\mathrm{T}}\bolds{\lambda}_0$,
where $\widehat{\bolds{\lambda}}(\bolds{\beta}_0, {\mathbf{w}}_0)$ solves (\ref{eqnstep1}), and\vspace*{-2pt}
define
\begin{eqnarray*}
{\sigma}^2(u, \mathbf{w}_0) &\equiv& \frac{1}{n}
\mathbf{B}_r(u)^{\mathrm{T}}E\bigl( \bigl[\widetilde {
\mathbf{Q}}_{\lambda i}\bigl\{ \mathbf{w} _0(\mathbf{T}_i)
\bigr\} ^{\mathrm{T}}\bolds{\Theta}_i\bigl\{\bolds{\beta
}_0, m_0, \mathbf{w}_0(\mathbf{T}_i)
\bigr\} \bolds{\Omega}^{-1}_i
\\[-2pt]
&& \hspace*{68pt}{}\times\bolds{\Theta}_i\bigl\{\bolds{\beta}_0,
m_0, \mathbf{w}_0(\mathbf{T}_i)\bigr\}
\widetilde{\mathbf{Q}}_{\lambda i}\bigl\{ \mathbf{w}_0(
\mathbf{T}_i)\bigr\} \bigr]\bigr) ^{-1 }
\\[-2pt]
&&{}\times E\bigl(\bigl[
\widetilde{\mathbf{Q}}_{\lambda i}\bigl\{\mathbf{w}_0(
\mathbf{T}_i)\bigr\} ^{\mathrm{T}}\bolds{\Theta}_i\bigl\{\bolds{\beta}_0,
m_0, \mathbf{w}_0(\mathbf{T}_i)\bigr\}
\bolds{\Omega}^{-1}_i \bolds{\Omega}_i^*
\bolds{\Omega}^{-1}_i\\[-2pt]
&&\qquad\hspace*{52pt}{}\times \bolds{\Theta}_i\bigl\{
\bolds{\beta}_0, m_0, \mathbf{w}_0(
\mathbf{T}_i)\bigr\}\widetilde{\mathbf{Q}}_{\lambda i}\bigl\{\mathbf{w}_0(
\mathbf{T}_i)\bigr\} \bigr]\bigr)
\\[-2pt]
&&{}\times E\bigl(\bigl[\widetilde{
\mathbf{Q}}_{\lambda i}\bigl\{ \mathbf{w} _0(\mathbf{T}_i)
\bigr\} ^{\mathrm{T}}\bolds{\Theta}_i\bigl\{\bolds{\beta
}_0, m_0, \mathbf{w}_0(\mathbf{T}_i)
\bigr\}
\\[-2pt]
&&\hspace*{27pt}{}\times \bolds{\Omega}^{-1}_i \bolds{
\Theta}_i\bigl\{\bolds{\beta}_0, m_0,
\mathbf{w}_0(\mathbf{T}_i)\bigr\} \widetilde{
\mathbf{Q}}_{\lambda i}\bigl\{ \mathbf{w}_0(\mathbf{T}_i)
\bigr\} \bigr]\bigr)^{-1 } \mathbf{B}_r(u),
\end{eqnarray*}
where $\bolds{\Omega}_i^* = E\{(\mathbf{D}_i - \mathbf{H}_i)^{\otimes2} |\mathbf{X}_i, \mathbf{Z}_i\}$ is the
true covariance matrix,\vspace*{-2pt} and
\begin{eqnarray*}
\sigma^2_w &\equiv& \frac{1}{n}
\mathbf{B}_r^{\mathrm{T}}(u) E \Biggl\{ \Biggl(
\mathbf{V}^{-1} E \Biggl[\sum_{k=1}^{M_i}
\sum_{v=1}^{M_i} E \bigl\{C_{ikv}
\Theta_{ik}\bigl\{\bolds{\beta}_0, m_0,
\mathbf{w}_0(T_{ik})\bigr\}
\\[-3pt]
&&\hspace*{135pt}{}\times\Theta_{iv}\bigl\{\bolds{\beta}_0,
m_0, \mathbf{w}_0(T_{iv})\bigr\}
\mathbf{B}_{r}\bigl\{\mathbf{w} _0(T_{iv})^{\mathrm{T}}
\mathbf{Z}_{iv}\bigr\}
\\[-2pt]
&&\hspace*{135pt}{}\times m_0'\bigl(
\mathbf{w}_0(T_{ik})^{\mathrm{T}}\mathbf{Z}_{ik}
\bigr)\\[-2pt]
&&\hspace*{135pt}{}\times\mathbf{Z}_{ik}^{\mathrm{T}} \bigl( \bigl\{\mathbf{B}(T_{ik})f_T(T_{ik})
\bigr\}^{-1} \\
&&\hspace*{42pt}\hspace*{123pt}{}\times\bigl[ {\mathbf{A}}_{wj}\bigl\{\bolds{
\beta}_0, m_0,\mathbf{w}_0(T_{ik})
\bigr\}\\[-2pt]
&&\hspace*{57pt}\hspace*{123pt}{}\times {\mathbf{V}}_{wj}\bigl\{\bolds{\beta}_0,
m_0,
\mathbf{w}_0(T_{ik})\bigr\}^{-1}\\
&&\hspace*{57pt}\hspace*{123pt}{}\times
\mathbf{K}_h(\mathbf{T}_j - T_{ik}) \\
&&\hspace*{57pt}\hspace*{136pt}{}\times{
\mathbf{S}}_{wj}\bigl\{\bolds{\beta}_0, m_0,
\mathbf{w} _0(T_{ik})\bigr\} \bigr] \bigr)
|
\\[-2pt]
&&\hspace*{40pt} \hspace*{57pt}\hspace*{123pt}{}M_i,\mathbf{O}_j \bigr\} |\mathbf{O}_j \Biggr]
\Biggr)^{\otimes2} \Biggr\} \mathbf{B}_r (u).
\end{eqnarray*}
Here
$\mathbf{V}$ is as defined in the notation of \textup{step}~1 in Appendix~\ref{secnotationEST},
and $C_{ikv} $ is the $(k, v)$th entry of the matrix $\bolds{\Omega}_i^{-1}$.
Then we\vspace*{-3pt} have
\begin{eqnarray*}
\bigl\{{\sigma}^2(u, \mathbf{w}_0) +
\sigma^2_w\bigr\}^{-1/2} \bigl(\widehat{m}
\bigl[u, \widehat{\bolds{\lambda}}\bigl\{\widehat{\bolds{\beta}}, {\widehat {
\mathbf{w}}}(\widehat{\bolds{\beta}})\bigr\}\bigr] - m_0(u) \bigr) &
\stackrel{d} {\rightarrow}& N(0, 1).
\end{eqnarray*}
Further because the order of $\sigma^2$ and $\sigma_\lambda^2$
are both $(nh_b)^{-1}$, together with Fact~1 in Section \textup{S.2}, we\vspace*{-2pt} have
\begin{eqnarray*}
\bigl|\widehat{m}\bigl[u, \widehat{\bolds{\lambda}}\bigl\{\widehat{\bolds{\beta}}, {
\widehat {\mathbf{w}}}(\widehat{\bolds{\beta}})\bigr\}\bigr]- m_0(u)\bigr|
&=& O_p\bigl\{ (nh_b)^{-1/2} +
h_b^q\bigr\},
\\[-2pt]
\bigl|\widehat{m}'\bigl[u, \widehat{\bolds{\lambda}}\bigl\{\widehat{
\bolds{\beta }}, {\widehat{\mathbf{w}}}(\widehat{\bolds{\beta}})\bigr\}\bigr] -
m_0'(u)\bigr| &=& O_p\bigl\{n^{-1/2}h_b^{-3/2}
+ h_b^{q-1}\bigr\}
\end{eqnarray*}
uniformly for\vadjust{\goodbreak} $ u \in(0, 1)$.
\end{Th}

Theorem~\ref{th3} shows that the estimation error of $\widehat{m}[u,
\widehat{\bolds{\lambda}}\{\widehat{\bolds{\beta}}, {\widehat
{\mathbf{w}}}(\widehat{\bolds{\beta}})\}]$ consists of two
components,
the approximation error of $\widehat{m}[u,
\widehat{\bolds{\lambda}}\{\widehat{\bolds{\beta}}, {\widehat
{\mathbf{w}}}(\widehat{\bolds{\beta}})\}]$ and the
approximation error of $\widetilde{m}(u,
\bolds{\lambda}_0)$, from their respective true functions. The errors of
$\widehat{m}$ and $\widehat{m}'$ go to
zero with the rates of $O_p\{(nh_b)^{-1/2}\}$ and
$O_p(n^{-1/2}h_b^{-3/2})$, respectively.
Under condition (A5), $\widehat{m}$ and
$\widehat{m}'$ are both consistent, and they approach the truths with the
standard B-spline convergence
rate.
We provide an outline of the proofs for Theorems \ref{th1}--\ref{th3}
in the supplementary article [\citep{jiang2015}]. The proofs are highly
technical and lengthy, and they require
several preliminary results, which we summarize as lemmas. We present
and prove these lemmas in the supplementary article [\citep{jiang2015}].

\section{Numeric evaluation via simulations}\label{secsim}

We now evaluate the finite sample performance of the proposed
estimation procedure on simulated data
sets. We simulate 1000 data sets from model (\ref{eqmodel})
under three settings. In Settings 1 and~2, we consider binary response and
use logit link function for $H$, while
in Setting 3, we consider continuous normal response and use an
identity $H$ function. In Setting 1, we
choose $m$ as a polynomial function with degree two. We
generate $\mathbf{w}$ initially as positive linear functions on $t$,
and then normalize
the vector to have summation one. Note that
the normalization function modifies the structure of $\mathbf{w}(t)$
and results
in a nonlinear vector-valued function in $t$. Additionally, we generate
$\mathbf{Z}_{ik}$
from the Poisson distribution and normalize the vectors by the sample
standard deviations. Furthermore, we generate $T_{ik}$ from the exponential
distribution and the covariate $\mathbf{X}_{ik}$ from the univariate normal
distribution.
In Settings 2 and 3, we use the sine
function for $m$, and generate $\mathbf{w}$ as power functions on $t$
and then
normalize the vector to have summation one. We generate
covariate vector $\mathbf{X}_i$ from a three-dimensional multivariate normal
distribution.
In order to stabilize the computation and control
numerical errors, in both settings, we transform the function $\mathbf{w}
(T_{ik})^{\mathrm{T}}\mathbf{Z}_{ik}$ to
$F\{\mathbf{w}(T_{ik})^{\mathrm{T}}\mathbf{Z}_{ik}\} = {\bolds\Phi} ( [\mathbf{w}
(T_{ik})^{\mathrm{T}}\mathbf{Z}_{ik} - E\{\mathbf
{w}^0(T_{ik})^{\mathrm{T}}\mathbf{Z}_{ik}\} ]
/\sqrt{\operatorname{var}\{\mathbf{w}^0(T_{ik})^{\mathrm{T}}\mathbf
{Z}_{ik}\}} )$, where $\mathbf{w}^0$
is the
initial value of $\mathbf{w}$, and
$E\{\mathbf{w}^0(T_{ik})^{\mathrm{T}}\mathbf{Z}_{ik}\}$ and
$\operatorname{var}\{\mathbf{w}^0(T_{ik})^{\mathrm{T}}\mathbf{Z}
_{ik}\}$
are approximated by the sample mean and the sample variance.
We then use B-spline to approximate $m\circ F^{-1}$ instead of $m$,
where $\circ$ denotes composite.
All other operations remain the same,
and the estimation and inference of the functional single index risk
score $m\{\mathbf{w}(T_{ik})^{\mathrm{T}}\mathbf{Z}_{ik}\}$, our
main research interest, is
carried out as described before.
To recover information regarding $m$, one can use the Delta method to obtain
the estimate and the variance of estimating
$m$ from that of estimating $m\circ F^{-1}$.

In all the implementations, we use the third
order quadratic spline. We select the number of internal knots $N =
\{n^ {1/5}(\log n)^{2} / 5\}$, which satisfies condition (A5) in Section~\ref{secthm}. We choose the Gaussian kernel with bandwidth $h =
n^{-2/15} h_s$, where $h_s$ is Silverman's rule-of-thumb bandwidth
[\citep{Silverman1986}]. Because $h_s = O(n^{-1/5})$, the
bandwidth selection satisfies condition (A2) in Section~\ref{secthm}.

Table~\ref{tabsim} shows the averaged point
estimators of $\bolds{\beta}$,
the empirical standard deviations calculated from the sample
variances, the averages of the estimated asymptotic standard deviation
($\Sigma^{1/2}$ in Theorem~\ref{th2}) over the
simulated samples and the mean squared errors ($\mbox{MSE}$) when the
sample sizes are 100, 500, 800, respectively. The conclusions are
similar under the three settings. To sum up, the
estimation biases are consistently small across all
samples sizes, the
empirical standard deviations and the estimated asymptotic
standard deviations are decreasing when the sample size increases. The
$\mbox{MSE}$
decreases as the sample size increases as well, mainly due to
the declining variations. Further, the empirical standard deviation of
the estimators
and average of the estimated standard deviations calculated from the
asymptotic results are close. In addition, the coverage
probabilities of the empirical confidence intervals are close to the
normal level of 95\%. This suggests that we can use the
asymptotic properties to perform inference and can obtain
sufficiently reliable results under moderate sample sizes.
%
\begin{sidewaystable}
\tabcolsep=0pt
\tablewidth=\textwidth
\caption{Simulation\vspace*{1pt} results in Settings 1, 2 and 3, based
on 1000 data sets. The true
parameter $\beta_0$, mean ($E$), empirical standard deviation [$\operatorname{sd}(\widehat\beta)$] and
average of the estimated standard deviations
[$\widehat{\operatorname{sd}}(\widehat\beta)$]
$\mbox{MSE}=\{{\operatorname{sd}}(\widehat{\beta}) \}^2+ \{E(\widehat
{\beta}) -
\beta\}^2$, the coverage probabilities (CP) of the 95\% empirical
confidence intervals are~reported}\label{tabsim}
\begin{tabular*}{\tablewidth}{@{\extracolsep{\fill}}lcccccccccccccccccc@{}}
\hline
& \multicolumn{6}{c}{\textbf{Setting 1}} & \multicolumn{6}{c}{\textbf{Setting 2}} &   \multicolumn{6}{c@{}}{\textbf{Setting 3}} \\[-4pt]
& \multicolumn{6}{c}{\hrulefill}  &  \multicolumn{6}{c}{\hrulefill} & \multicolumn{6}{c}{\hrulefill}\\
& $\bolds{\beta_0}$& $\bolds{E(\widehat{\beta})}$& $\bolds{{\operatorname{sd}}(\widehat{\beta
})}$&$\bolds{\widehat{\operatorname{sd}}(\widehat{\beta})}$& \textbf{MSE} &\textbf{CP} &
$\bolds{\beta_0}$&$\bolds{E(\widehat{\bolds{\beta
}})}$&$\bolds{{\operatorname{sd}}(\widehat{\bolds{\beta}})}$&$\bolds{\widehat
{\operatorname{sd}}(\widehat{\bolds{\beta}})}$&\textbf{MSE}&\textbf{CP}& $\bolds{\beta_0}$ & $\bolds{E(\widehat{\bolds{\beta
}})}$&$\bolds{{\operatorname{sd}}(\widehat{\bolds{\beta}})}$&$\bolds{\widehat{\operatorname{sd}}(\widehat{\bolds{\beta}})}$&\textbf{MSE}&\textbf{CP}
\\
\hline
& \multicolumn{18}{c@{}}{$n = 100$}\\
$\beta_1$&$-$0.2&$-$0.202&0.157&0.113&0.0247&0.957  &$-$0.5&$-$0.505&0.131&0.116&0.0171&0.908 &$-$0.5&$-$0.501&0.062&0.052&3.85\mbox{e--}3&0.938 \\
$\beta_2$& $-$0.4&$-$0.398&0.119&0.115&0.0142&0.940 &\phantom{$-$}0.2&$-$0.200&0.122&0.112&0.0147&0.923  &\phantom{$-$}0.2&$-$0.200&0.060&0.053&3.60\mbox{e--}3&0.932  \\
$\beta_3$&$-$0.6&$-$0.601&0.124&0.118&0.0153&0.957  &\phantom{$-$}0.5&$-$0.515&0.125&0.116&0.0159&0.927   &\phantom{$-$}0.5&$-$0.503&0.061&0.053&3.73\mbox{e--}3&0.932 \\[3pt]
& \multicolumn{18}{c@{}}{$n = 500$}\\
$\beta_1$&$-$0.2&$-$0.198&0.052&0.050&0.0027&0.954& $-$0.5&$-$0.507&0.056&0.053&0.0032&0.946&$-$0.5&$-$0.500&0.025&0.024&6.25\mbox{e--}4&0.966 \\
$\beta_2$&$-$0.4&$-$0.398&0.053&0.051&0.0028&0.947&\phantom{$-$}0.2&$-$0.198&0.053&0.052&0.0028&0.951  &\phantom{$-$}0.2&$-$0.200&0.024&0.024&5.76\mbox{e--}4&0.945  \\
$\beta_3$&$-$0.6&$-$0.601&0.056&0.053&0.0031&0.939 &\phantom{$-$}0.5&$-$0.508&0.054&0.053&0.0031&0.944 &\phantom{$-$}0.5&$-$0.502&0.025&0.024&6.29\mbox{e--}4&0.963  \\[3pt]
& \multicolumn{18}{c@{}}{$n = 800$}\\
$\beta_1$&$-$0.2&$-$0.197&0.041&0.040&0.0017&0.951 &  $-$0.5&$-$0.504&0.043&0.042&0.0019&0.953& $-$0.5&$-$0.500&0.020&0.019&4.00\mbox{e--}4&0.949 \\
$\beta_2$&$-$0.4&$-$0.398&0.041&0.040&0.0017&0.949  & \phantom{$-$}0.2&$-$0.202&0.041&0.041&0.0017&0.962  &\phantom{$-$}0.2&$-$0.200&0.019&0.019&3.61\mbox{e--}4&0.949 \\
$\beta_3$&$-$0.6&$-$0.602&0.044&0.042&0.0019 &0.946 & \phantom{$-$}0.5&$-$0.505&0.043&0.042&0.0019&0.951 &\phantom{$-$}0.5&$-$0.501&0.020&0.019&4.01\mbox{e--}4&0.952   \\
\hline
\end{tabular*}
\end{sidewaystable}

\begin{figure}[p]

\includegraphics{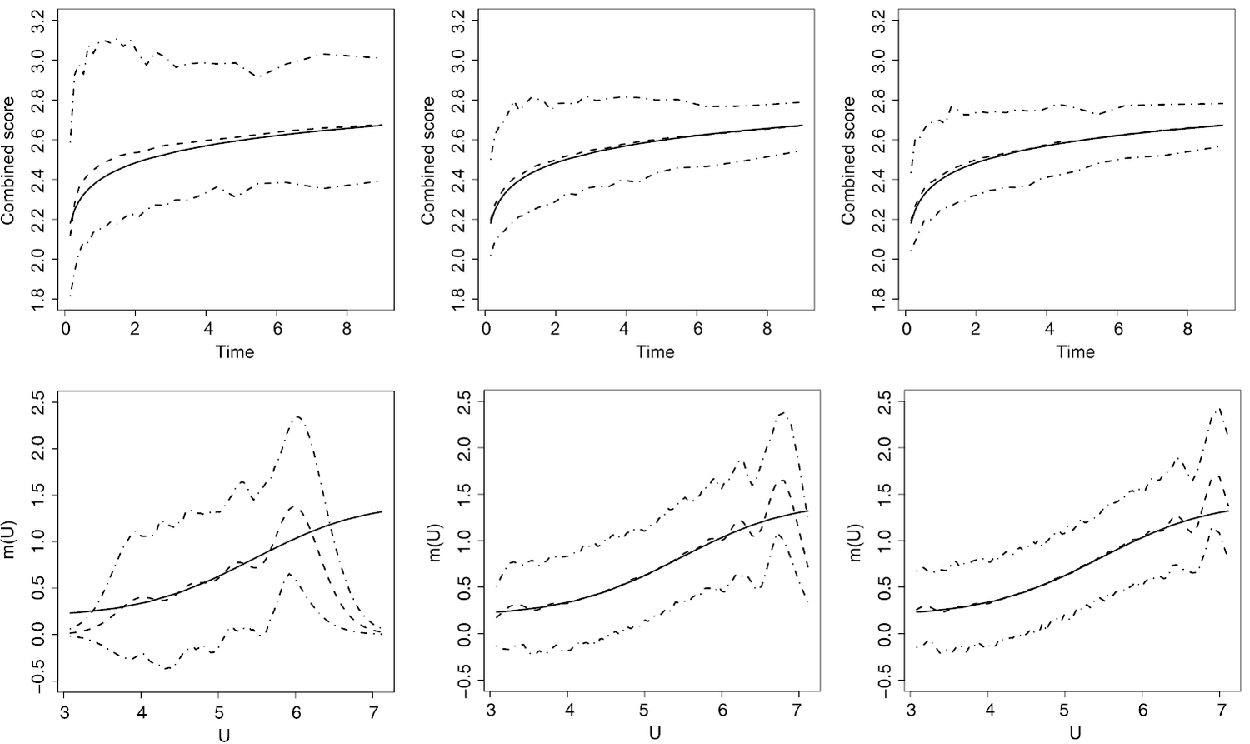}

\caption{Estimation of ${\mathbf{w}}(t)^{\mathrm{T}}\mathbf{z}$ (upper)
and $m(u)$
(bottom) as a function of $t$ and $u$, respectively, in
Setting 1 with
sample sizes 100 (left), 500 (middle) and 800 (right).
True function (solid line), average of 1000 estimated functions
(dashed lines) and 95\% pointwise confidence band (dash-doted lines)
are provided.} \label{figsimuw}
\end{figure}

\begin{figure}[p]

\includegraphics{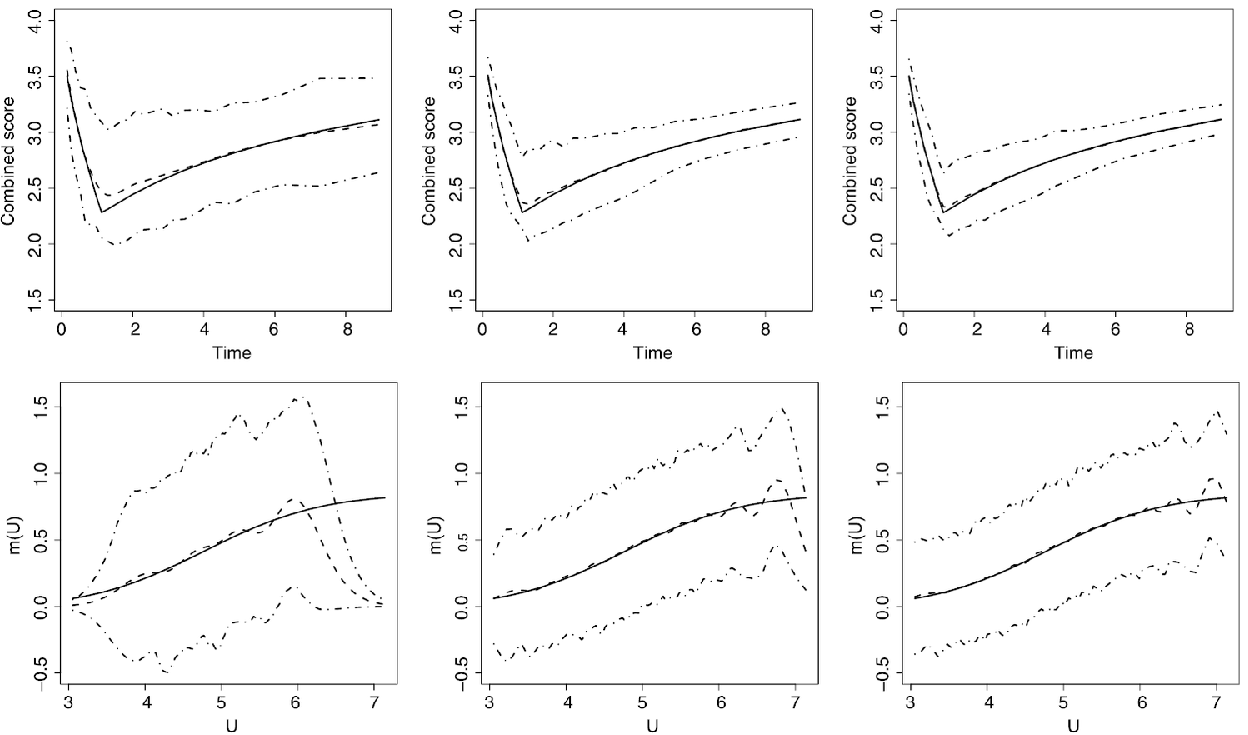}

\caption{Estimation of ${\mathbf{w}}(t)^{\mathrm{T}}\mathbf{z}$ (upper) and
$m(u)$ (bottom) as a function of $t$ and $u$, respectively,
in Setting 2 with
sample sizes 100 (left), 500 (middle) and 800 (right).
True function (solid line), average of 1000 estimated functions
(dashed lines), and 95\% pointwise confidence band (dash-doted lines)
are provided.} \label{figsimuwsin}
\end{figure}

\begin{figure}[t]

\includegraphics{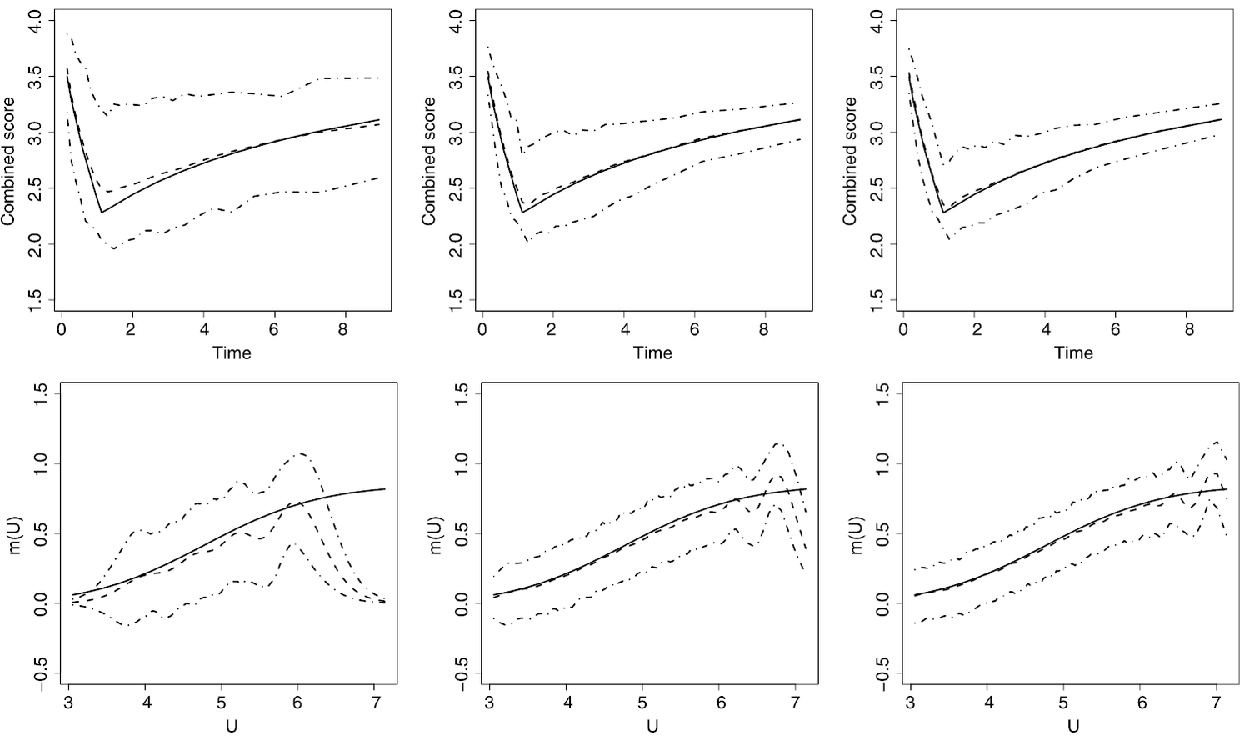}

\caption{Estimation of ${\mathbf{w}}(t)^{\mathrm
{T}}\mathbf{z}$ (upper) and
${m}(u)$ (bottom) as a function of $t$ and $u$, respectively,
in Setting 3 with
sample sizes 100 (left), 500 (middle) and 800 (right).
True function (solid line), average of 1000 estimated functions
(dashed lines) and 95\% pointwise confidence band (dash-doted lines)
are provided.} \label{figsimuwsinnorm}
\end{figure}

We also examined the performances of $\widehat{\mathbf{w}}$ and
$\widehat{m}$ to assess the properties of the estimated functional
single index risk score.
Under the first setting, because the functional single index risk
score is fixed with respect to $\beta$, we only evaluate the settings
with $\beta=
-0.4$.
To evaluate the combined score $\widehat{\mathbf{w}}(t) ^{\mathrm
{T}}\mathbf{Z}$ as
a function of~$t$, we fix $\mathbf{Z}$ at $\mathbf{Z}^*= (1, 2, 3,
4)$ and plot the
averages of the estimated combined score
$\widehat{\mathbf{w}}(t) ^{\mathrm{T}}\mathbf{Z}^*$ over the 1000
simulations
around the true scores ${\mathbf{w}}_0(t)
^{\mathrm{T}}\mathbf{Z}^*$ in the upper panels of Figures~\ref{figsimuw}, \ref
{figsimuwsin} and \ref{figsimuwsinnorm} for
Settings 1, 2 and 3, respectively. Additionally, we present the 95\%
pointwise confidence band. The results show
that the estimates are
close to the true function. Further, the 95\%
confidence band becomes narrower when the sample size increases, which
indicates that the estimation variation decreases with increased sample
size. Moreover, we evaluated the coverage probabilities of the
empirical pointwise confidence bands of $\mathbf{w}$, by computing the coverage
probabilities at a set of fixed points across $t$ and taking their average.
The average coverage probabilities for $n= 100, 500,
800$ are $0.934, 0.936, 0.939$ in Setting 1, $0.939, 0.940, 0.941$
in Setting 2 and $0.931, 0.934, 0.936 $ in Setting 3, respectively.
All are reasonably close to the
nominal level of 95\%.

To evaluate the performance of $\widehat{m}$, we plot the average of
$\widehat{m}(u)$ based on the 1000 simulations, as well as the 95\%
pointwise confidence band in the bottom panels of Figures~\ref{figsimuw},
\ref{figsimuwsin} and \ref{figsimuwsinnorm} for Settings 1, 2 and
3, respectively.
The plots show that the estimators are close to the
true functions except on the boundary when the sample size is relatively
small. In addition, when the sample size increases, the confidence band
becomes narrower, benefiting from the smaller estimation
variation. Note that because of the additional transformation on
$\mathbf{w}(t)^{\mathrm{T}}\mathbf{Z}$, it is expected that the true
$m$ function
does not appear to be periodic sine function on $\mathbf
{w}(t)^{\mathrm{T}}\mathbf{Z}$.
Moreover, we evaluate the converge probability of the empirical
pointwise confidence bands of $m$. The average coverage probabilities
are $0.943, 0.947, 0.948$ in Setting 1,
$0.957, 0.960, 0.951$ in Setting 2 and $0.939, 0.947, 0.946$ in
Setting 3, respectively. Again, they are all fairly close to the
nominal level of 95\%.

In summary, Table~\ref{tabsim}, Figures~\ref{figsimuw}, \ref{figsimuwsin}, \ref{figsimuwsinnorm} illustrate
the desirable finite sample performance of the fused kernel/B-spline
combination method in estimating
$\bolds{\beta}, m$ and~$\mathbf{w}$. In terms of parameter
estimation and
function estimation in the nonboundary region, the estimators show
very small biases
across all sample sizes, and
decreasing variability as the sample size increases.
The asymptotic variance and sample
empirical variance in estimating $\bolds{\beta}$ are close.
Furthermore, the coverage
probability of the empirical confidence intervals for $\bolds{\beta}$
and the
coverage probability of the empirical pointwise confidence bands for
$\mathbf{w}$ and $m$ are close to the
nominal levels, which supports using the asymptotic
results for the subsequent inferences.

\section{Application}\label{secreal}

We apply the functional single index risk score model and the fused
kernel/B-spline
semiparametric
estimation method to analyze a real data set from a Huntington's
disease (HD)
study.
Current research in HD aims to find reliable
prodromes to enable early detection of HD. The joint effect of the
cognitive scores on odds of HD diagnosis is shown to
change with time. In addition,
the relationship between the cognitive symptoms and the log-odds of the
disease diagnosis is shown to be nonlinear [\citep{Paulsen2008}].
Our goal is to study the nonlinear time dependent
cognitive effects so as to facilitate the early detection of~HD.

Specifically, let $D_{ik}$, $\mathbf{Z}_{ik}$ and $\mathbf{X}_{ik}$
represent the
binary disease
indicator, the cognitive
score vector and the additional covariate vector for the $i$th
individual at the $j$th
measurement time,
respectively. The cognitive scores include
SDMT [\citep{Smith1982}], Stroop color, Stroop word and Stroop
interference tests [\citep{Stroop1935}]. They are denoted by $Z_{i1},
\ldots, Z_{i4}$,
respectively. The covariates of interest are gender,
education, CAP score [\citep{Zhang2011}].
They
are denoted by $X_{i1}, \ldots, X_{i3}$, respectively. The subject's
age at the visiting time
serves as the time variable $T_{ik}$. We normalize the continuous
variables to the interval $(0, 1)$ to alleviate numerical
instability. Without changing notation, we transform $Z_{i1}, \ldots,
Z_{i4}$, $X_{i3}$, $T_{ik}$ by the normal distribution functions with
means and variances estimated from the sample.

We use logit link function to model the binary outcomes; that is,
we assume
%
\begin{eqnarray}
\label{eqnlgtmdl}
&& H\bigl[m\bigl\{\mathbf{w}(T_{ik})^{\mathrm{T}}
\mathbf{Z}_{ik}\bigr\} +\bolds{\beta}^{\mathrm{T}}\mathbf{X}_{ik}
\bigr] = \frac{\exp
 [m\{\mathbf{w}(T_{ik})^{\mathrm{T}}\mathbf{Z}_{ik}\}+\bolds
{\beta}^{\mathrm{T}}\mathbf{X}_{ik} ]}{1 +
\exp
 [m\{\mathbf{w}(T_{ik})^{\mathrm{T}}\mathbf{Z}_{ik}\}+\bolds
{\beta}^{\mathrm{T}}\mathbf{X}_{ik} ]}.
\end{eqnarray}
We obtain the initial estimates and a working correlation matrix
using the GEE method with exchangeable
covariance assumption. We
choose the exchangeable covariance structure because in our
setting, it facilitates computation
and accounts for the longitudinal correlations. Let the working
correlation coefficient matrix be $\mathbf{R}_i$, the working covariance
matrix be
$\widehat{\Theta}_i ^ {1/2} \mathbf{R}_i \widehat{\Theta}_i^
{1/2}$, where
$\widehat{\Theta}_i $ is $\mathbf{H}_i (1 - \mathbf{H}_i)$ with
estimated $\widehat{\bolds{\lambda}
}, \widehat{\mathbf{w}}, \widehat{\bolds{\beta}}$
plugged in. We implement the profiling
procedure described in Section~\ref{secprofile} in the subsequent
estimation. The kernel and B-spline functions are defined in the
same way as described in Section~\ref{secsim}. We obtain the point\vspace*{1pt}
estimators $\widehat{\bolds{\beta}} = (-0.34, -0.89,
2.31)^{\mathrm{T}}$ and the asymptotic variances $\operatorname
{var}(\widehat{\bolds{\beta}})= (0.0035,
0.00044,0.011)^{\mathrm{T}}$. Consequently, the
95\% asymptotic confidence intervals are $\{(-0.46, -0.23)$,
$(-0.93, -0.85)$, $(2.09, 2.52)\}$, which demonstrate
the significant effect of gender, education level and CAP score on
the disease risk.
Specifically, 
females ($X_{i1} =0$) tend to have higher disease risk than males ($X_{i1}
= 1$). In addition, patients with lower education levels and higher CAP
scores are more likely
to develop Huntington's disease, which is consistent with the clinical
literature [\citep{Zhang2011}].

\begin{figure}

\includegraphics{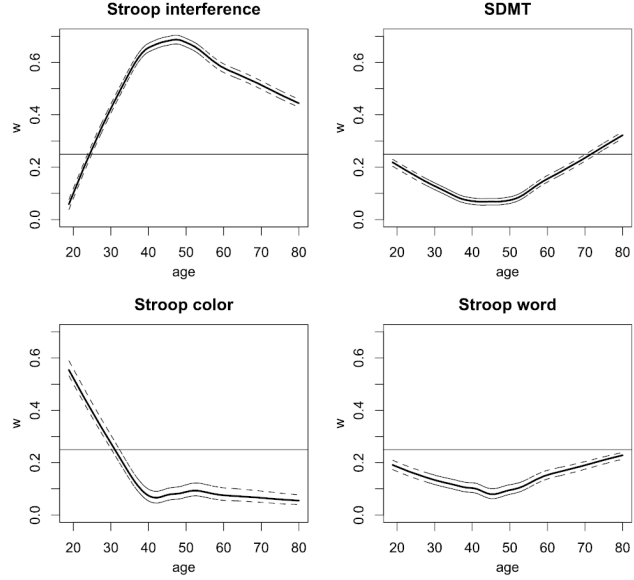}

\caption{Estimation of the weight function $w(t)$'s and
the $95\%$ asymptotic confidence bands in Huntington's disease data.
The reference line is $0.25$.}
\label{figwplot}
\end{figure}

We also plot $\widehat{\mathbf{w}}(t)$ to show the variation
patterns of the effect of the four cognitive scores over time. Figure~\ref{figwplot} shows that the Stroop interference
score has a more important effect than all the others after age 30. The
95\% pointwise
confidence interval remains above the 0.25 level after age 27,
and the Stroop interference
score effect largely dominates all the other effects during that
period. This dominating effect
indicates that the Stroop inference
score has the closest relationship with the onset of HD, and in turn could
be used to predict HD most effectively among the
four. Further, Stroop color has a large effect at earlier ages (before
30 or
at early 30's), while the
SDMT has a reasonably large effect at later ages (75 or above).
Moreover, Stroop word
has relatively small predicative
effects (${<}0.25$) on the disease risk across all ages.
The plots clearly show the time dependent
nature of the cognitive score effects. More
specifically, Stroop color effect is decreasing over times, Stroop
interference effect is a concave function of
time, while SDMT, Stroop word effects are convex functions of
time. The last three nonmonotone effects reach their extreme values
around the ages of 40 to 50.
In summary, the results show that the
Stroop interference is more
relevant to the disease risk than the other scores.
Further, the relative magnitude of the score effects clearly change
over time, which suggests the need
to closely monitor specific cognitive scores for different age groups.
This illustrates the importance of modeling $\mathbf{w}$ as a function
of age, and the convenience of using a weighted score $\mathbf{w}(t)^{\mathrm{T}}\mathbf{Z}$
as a combined cognitive profile
in practice.

\begin{figure}

\includegraphics{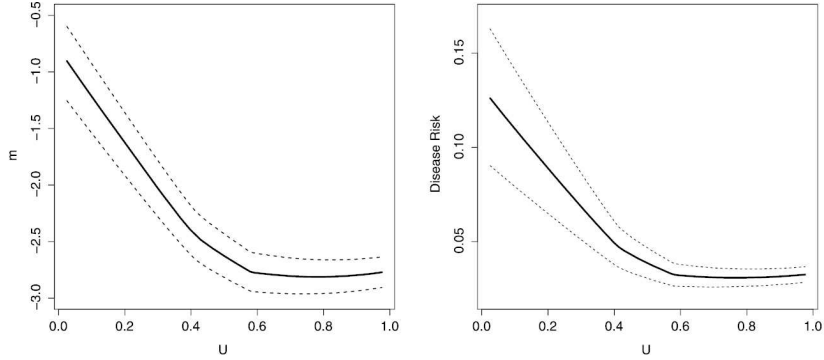}

\caption{Function $\widehat{m}(u)$ (left) and the estimated disease
risk as
a function of $u$ (right) in Huntington's disease
data.}
\label{figfudis}
\end{figure}

The form of the function
$\widehat m$ is shown in the left panel of Figure~\ref{figfudis}. We also plot the 95\% pointwise asymptotic confidence
band of $\widehat m$ in the range of the combined scores $U$. The plot shows
that the
functional single index risk score is a decreasing function of the
index. The upper confidence
interval does not include $0$, which shows that the functional single
index risk score
is significantly smaller than 0 at any age and cognitive score values
in this population.

In the right panel of Figure~\ref{figfudis}, we plot the disease risk
(the estimated
probability of $D=1$) and the 95\% pointwise asymptotic confidence
band, where the confidence band is based on estimated
variance, calculated using the Delta method and the estimated variance
of $\widehat m$.
The results show that the disease risk
decreases with the combined cognitive score value $U$. The 95\%
confidence interval does not include
the 0.5 line, which shows that the disease risk in the population is
smaller than 0.5 across all age and cognitive score values. Combining
the two plots, Figure~\ref{figfudis} shows that a higher value of the combined score
$U=\mathbf{w}(t)^{\mathrm{T}}\mathbf{Z}$, which implies better
cognitive functioning, tends
to lower functional single index risk score
and in turn lower the risk of HD. The effect of the functional single
index cognitive risk score on HD diagnosis is approximately quadratic
for a standardized score $U<0.6$, and is approximately a constant for
$U>0.6$. The flattening of the effect reflects a ceiling effect for
subjects with better cognitive performance.

\begin{figure}

\includegraphics{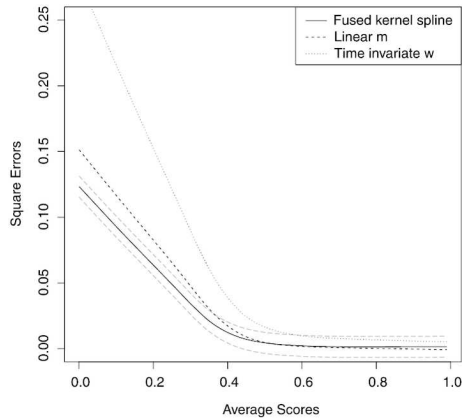}

\caption{The mean squared predictive errors versus the standardized
averaged score $\sum_{j= 1}^4
Z_{ik}/4$ in Huntington's disease data. The gray lines are the 95\%
confidence intervals for the fused
kernel/B-spline method.}\label{figprederrors}
\end{figure}

Next, we perform two sensitivity analyses to justify using a more flexible
generalized
partially linear functional single index model as shown in
(\ref{eqnlgtmdl}). We
compare model (\ref{eqnlgtmdl}) with two simpler models. The first
assumes the function $m$ is linear, hence
%
\begin{eqnarray}
\label{eqnmdl2}
&& H\bigl\{\mathbf{X}_{ik},\mathbf{Z}_{ik};{
\bolds\theta},\mathbf{w}(T_{ik})\bigr\}=\frac{\exp\{
\alpha_c+\alpha_1\mathbf{w}(T_{ik})^{\mathrm{T}}\mathbf{Z}_{ik}+\bolds{\beta}^{\mathrm{T}}\mathbf{X}_{ik}\}}{
1+\exp\{\alpha_c+\alpha_1\mathbf{w}(T_{ik})^{\mathrm{T}}\mathbf{Z}_{ik}+\bolds{\beta}^{\mathrm{T}}\mathbf{X}
_{ik}\}},
\end{eqnarray}
where $\alpha_c, \alpha_1$ are unknown parameters. The second
assumes the weight function $\mathbf{w}$ is time-invariant,
hence
%
\begin{eqnarray}
\label{eqnmdl3}
&&  H(\mathbf{X}_{ik},\mathbf{Z}_{ik};{
\bolds\theta},\mathbf{w})=\frac{\exp\{m(\mathbf{w}
^{\mathrm{T}}\mathbf{Z}_{ik})+\bolds{\beta}^{\mathrm{T}}\mathbf{X}_{ik}\}}{
1+\exp\{m(\mathbf{w}^{\mathrm{T}}\mathbf{Z}_{ik})+\bolds{\beta
}^{\mathrm{T}}\mathbf{X}_{ik}\}},
\end{eqnarray}
where $\mathbf{w}$ is an unknown parameter vector. We carried out the
estimation of
$\mathbf{w}(t)$ in the first model using kernel
method and the estimation for $m$ in the second model
via B-spline method. We implemented 1000 5-fold cross
validation analysis. We evaluated models by the mean squared predictive
error (i.e., the mean squared differences between $D_i$ and the
predicted probability of $D_i=1$ on the test set) as a function of the
average of the four standardized
cognitive scores $\sum_{j}^4 Z_j/4$, which we named the standardized
score. In Figure~\ref{figprederrors}, we plot the mean squared predictive error curves
obtained under the proposed model (\ref{eqnlgtmdl}) and two simpler models.
The results show that
our original generalized
partially linear model with functional single index
outperforms
model (\ref{eqnmdl3}) uniformly across the range of the standardized
scores in terms of a lower mean squared error. We also plot the
empirical 95\% confidence intervals
of the squared predictive errors under the proposed model. Compared
with the simpler model (\ref{eqnmdl2}), our model gives
significant smaller predictive errors when the standardized score is
smaller than 0.36. The medians of the squared
predicative errors in this range are 0.040 and 0.049 for models (\ref
{eqnlgtmdl}) and
(\ref{eqnmdl2}), respectively. When the standardized score is greater than
0.5, model (\ref{eqnmdl2}) performs slightly, but not significantly
better than model (\ref{eqnlgtmdl}). Overall, the total mean squared error
summarized by the area under the predictive error curves for models
(\ref{eqnlgtmdl}), (\ref{eqnmdl2}) and (\ref{eqnmdl3}) are, respectively,
0.022, 0.028 and 0.057, which
justifies using the more flexible model in
(\ref{eqnlgtmdl}) to fit the
Huntington's disease data. The results also demonstrate the potential
of using our method as an exploratory tool to assess general patterns
of data.

\section{Conclusion and discussions}\label{seccon}

We have developed a generalized partially linear functional
single index risk score model in the longitudinal data
framework. We
explore the relationship between the cognitive scores and the disease
risk so as to predict HD diagnosis early, and in turn to intervene with
the disease
progression in a timely manner.

We introduce a framework of jointly using the B-spline and kernel methods
in semiparametric estimation. We use
B-spline to approximate the functional single index
risk score
function $m$, and use a kernel smoothing technique for estimating the
cognitive weight functions of time $\mathbf{w}(t)$. We
integrate B-spline basis expansion, kernel smoothing and longitudinal
analysis, and have proven the consistency and asymptotic
normalities of the covariate coefficient
estimators, the time dependent weight function estimators and the
single index risk score function estimators. The derivation relies on the
assumption that the iteration procedure converges to a parameter vector
value that is in a small neighborhood of the truth, which generally
requires the estimating equation to have a unique zero. The unique
zero property is difficult to guarantee in theory and is less likely to
hold when sample size is small or moderate. To this end, empirical
knowledge is usually used to select a suitable root. In our
simulations, multiple roots issues did not occur, and the
numerical results show desirable finite sample properties of the
estimators. The real data analysis yields results which are
interpretable and useful in practice. In
summary, the functional single index model provides rich and
meaningful information regarding the association between the disease
risk and the cognitive score profiles. It is of course
also possible to use B-spline or kernel methods to estimate both
$m$ and $\mathbf{w}(t)$, and research along this line can also be interesting.

Our method accommodates both continuous and categorical
response variables as
long as the link function $H$ is continuously
differentiable and has finite second derivative.
One outstanding research question in these models, even in the
context when the marginal model is completely parametric (e.g.,
both $m$ and $\mathbf{w}$ are known), is the estimation efficiency. As
far as we
are aware, there is no guarantee that GEE family contains the efficient
estimator, and how to obtain an asymptotically efficient estimator is
certainly worth further research.

The proposed generalized partially linear functional single index model can
be used to incorporate high-dimensional data,
since the single index risk score is a natural method of alleviating the
curse of the dimensionality. For example, the single index score could
be a combination of gene expression covariates to facilitate the
genetic association study. 
Furthermore, the generalized partially linear functional single index
risk score can be used
in an adaptive randomization clinical trial study to improve study
efficiency. For example, we can use a single index
risk score to summarize some disease
related biomarkers, which provide early information about the primary
endpoints in adaptive trials. When a trial progresses, the
information can be used to make certain intermediate decisions, such as
treatment assignments among the patients, and stopping or continuation
of the
trial.

\begin{appendix}
\section{Proof of Proposition~\texorpdfstring{\lowercase{\protect\ref{proidentify}}}{1}}\label{secidentify}

Assume there exist $m_1 \in\mathcal{M}, \mathbf{w}_1(t) \in\mathcal
{D} $
and $\bolds{\beta}_1\in
\mathbb{R}^{d_{\beta}}$, such that
%
\begin{eqnarray}
\label{eqiden}
&&
m_1\bigl\{\mathbf{w}_1^{\mathrm{T}}(t)
\mathbf{Z}\bigr\} +\bolds{\beta }_1^{\mathrm{T}}\mathbf{X}=
m_0\bigl\{\mathbf{w}_0^{\mathrm{T}}(t) \mathbf{Z} \bigr
\} +\bolds{\beta}_0^{\mathrm{T}}\mathbf{X},
\end{eqnarray}
where $m_0, \mathbf{w}_0(t)$ and $\bolds{\beta}_0$ are the true
parameter values.
Taking derivative with respect to $\mathbf{Z}$ and $t$ on both sides of
the equation, we
obtain
%
\begin{eqnarray}
 m'_1\bigl\{\mathbf{w}_1^{\mathrm{T}}(t)
\mathbf{Z}\bigr\} \mathbf{w}_1(t) &=& m'_0
\bigl\{\mathbf{w}_0^{\mathrm{T}}(t) \mathbf{Z}\bigr\} \mathbf{w}
_0(t),
\nonumber
\\[-8pt]
\label{eqnderm}
\\[-8pt]
\nonumber
m'_1\bigl\{\mathbf{w}_1^{\mathrm{T}}(t)
\mathbf{Z}\bigr\} \mathbf{w}_1'(t)^{\mathrm{T}}
\mathbf{Z}&=& m'_0\bigl\{\mathbf{w}_0^{\mathrm{T}}
(t) \mathbf{Z}\bigr\} \mathbf{w}_0'(t)^{\mathrm{T}}
\mathbf{Z}.
\end{eqnarray}
Because $m_1, m_0$ are one-to-one,
$m'_1\{\mathbf{w}_1^{\mathrm{T}}(t) \mathbf{Z}\}=m'_0\{\mathbf{w}_0^{\mathrm{T}}(t) \mathbf{Z}\}=0$ can hold only
for a set of discrete set of $\mathbf{w}_1^{\mathrm{T}}(t) \mathbf{Z}$ and $\mathbf{w}_0^{\mathrm{T}}(t)
\mathbf{Z}$
values, hence a discrete set of
$t$ values. Thus due to the continuity of
$m_1', m_0', \mathbf{w}_1$ and $\mathbf{w}_0$,
(\ref{eqnderm}) implies $ \mathbf{w}_1'(t)^{\mathrm{T}}\mathbf{Z}/w_{1j}(t) =
\mathbf{w}_0'(t)^{\mathrm{T}}\mathbf{Z}/w_{0j}(t)$
for all $j=1, \dots, d_w$, all $\mathbf{Z}$ and all $t\in
[0, \tau]$.
Thus
$\mathbf{w}_1'(t)^{\mathrm{T}}E(\mathbf{Z}^{\otimes2})/w_{1j}(t) =
\mathbf{w}_0'(t)^{\mathrm{T}}E(\mathbf{Z}^{\otimes2})/w_{0j}(t)$.
Furthermore, $E(\mathbf{Z}
^{\otimes2})$ is positive
definite, and in turn it is invertible and leads to
${\mathbf{w}'_{1}(t)}/{w_{1j}(t)} =
{\mathbf{w}'_{0}(t)}/{w_{0j}(t)}$.
In particular, we have
$w'_{1j}(t)/w_{1j}(t) = w'_{0j}(t)/w_{0j}(t)$
for all $j = 1, \ldots, d_w$. This gives
$ w_{1j}(t) = w_{0j}(t) c_j$
for some constant $c_j$, or equivalently,
$\mathbf{w}_1(t) = \mathbf{C}\mathbf{w}_0(t)$
where $ \mathbf{C}$ is a diagonal matrix with $ c_j$'s on the diagonal.
Taking derivative with respect to $ t$, we further have
$\mathbf{w}_1'(t)=\mathbf{C}\mathbf{w}'_0(t)$. Dividing $w_{1j}(t)$
on both sides, we have
$\mathbf{w}_1'(t)/w_{1j}(t)=(\mathbf{C}/c_j) \mathbf{w}'_0(t)/w_{0j}(t)$.
Therefore, $\mathbf{C}/c_j$ is the identity matrix. In other
words, $c_j, j=1, \dots, d_w$ are identical. Since $\|\mathbf{w}_1(t)\|_1=\|
\mathbf{w}_0(t)\|_1=1$
and $\mathbf{w}_1(t)$, $\mathbf{w}_0(t)$ are positive,
this further implies $\mathbf{w}_1(t)=\mathbf{w}_0(t)$.
Therefore, (\ref{eqnderm}) reduces to
$m_1'\{\mathbf{w}_0^{\mathrm{T}}(t)
\mathbf{Z}\} - m_0'\{\mathbf{w}_0^{\mathrm{T}}(t) \mathbf{Z}\} =0$.
This further implies $m_1\{\mathbf{w}_0^{\mathrm{T}}(t)
\mathbf{Z}\} = m_0\{\mathbf{w}_0^{\mathrm{T}}(t) \mathbf{Z}\} + C_1$
for a
constant $C_1$. Because $m_1(0) = m_0(0) = c_0$,
$C_1 = 0$, that is, $m_1= m_0$. Model
(\ref{eqiden}) now leads to
$\bolds{\beta}_1^{\mathrm{T}}\mathbf{X}= \bolds{\beta}_0^{\mathrm
{T}}\mathbf{X}$. The equality holds for any $\mathbf{X}$,
which implies $\bolds{\beta}_1^{\mathrm{T}}E(\mathbf{X}^{\otimes
2}) =
\bolds{\beta}_0^{\mathrm{T}}E(\mathbf{X}^{\otimes2})$. Since
$E(\mathbf{X}^{\otimes2})$ is
positive definite, and in turn is
invertible, we have $\bolds{\beta}_1= \bolds{\beta}_0$. Therefore,
we have $\bolds{\beta}_1 =
\bolds{\beta}_0$, $\mathbf{w}_1(t)=\mathbf{w}_0(t)$ and $m_1
= m_0$, and hence
the problem is identifiable.

\section{Notation in estimation step}\label{secnotationEST}

\subsection*{Notation in step 1}
We define an $M_i \times d_\lambda$ matrix
\begin{eqnarray*}
&&  \widetilde{\mathbf{Q}}_{\lambda i}\bigl\{\mathbf{w}(\mathbf{T}_i)
\bigr\} = \lleft[
\matrix{
B_{r1}\bigl\{
\mathbf{w}(T_{i1})^{\mathrm{T}}\mathbf{Z}_{i1}\bigr\}& \ldots
&B_{rd_\lambda}\bigl\{\mathbf{w} (T_{i1})^{\mathrm{T}}
\mathbf{Z}_{i1}\bigr\}
\cr
\vdots&\vdots&\vdots
\cr
B_{r1}\bigl\{\mathbf{w}(T_{iM_i})^{\mathrm{T}}
\mathbf{Z}_{iM_i}\bigr\}& \ldots &B_{rd_\lambda}\bigl\{\mathbf{w}
(T_{iM_i})^{\mathrm{T}} \mathbf{Z}_{iM_i}\bigr\}}
 \rright],
\end{eqnarray*}
and define $\widetilde{\mathbf{Q}}_{\lambda i}\{\mathbf{w}(t_0)\}$
to be the same as $
\widetilde{\mathbf{Q}}_{\lambda i}\{\mathbf{w}(\mathbf{T}_i)\} $,
except we replace $T_{ik}, k = 1,
\ldots, M_i$ with $t_0$.
Here and throughout the text, replacing $\mathbf{T}_i$ by $t_0$ means
that we replace
$T_{ik}= t_0$ for each $k, k = , 1, \ldots, M_i$.
Let
\begin{eqnarray*}
\mathbf{V}_n &=& n^{-1} \sum
_{i=1}^{n}\bigl[\widetilde{ \mathbf
{Q}}_{\lambda i}\bigl\{ \mathbf{w}_0(\mathbf{T}_i)
\bigr\} ^{\mathrm{T}} \bolds{\Theta}_i \bigl\{\bolds{
\beta}_0, m_0, \mathbf{w}_0(
\mathbf{T}_i)\bigr\} \bolds{\Omega }^{-1}_i
\\
&&\hspace*{48pt}{}\times\bolds{\Theta}_i \bigl\{\bolds{\beta}_0,
m_0, \mathbf{w}_0(\mathbf{T}_i)\bigr\}
\widetilde {\mathbf{Q}}_{\lambda i}\bigl\{ \mathbf{w}_0(
\mathbf{T}_i)\bigr\} \bigr],
\\
\mathbf{V}&=&E \bigl( \bigl[\widetilde{\mathbf{Q}}_{\lambda i}\bigl\{ \mathbf
{w}_0(\mathbf{T}_i)\bigr\} ^{\mathrm{T}}\bolds{
\Theta}_i \bigl\{\bolds{\beta}_0, m_0,
\mathbf{w}_0(\mathbf{T}_i)\bigr\} \bolds{\Omega
}^{-1}_i
\\
&&\hspace*{29pt}{}\times\bolds{\Theta}_i \bigl\{\bolds{\beta}_0,
m_0, \mathbf{w}_0(\mathbf{T}_i)\bigr\}
\widetilde {\mathbf{Q}}_{\lambda i}\bigl\{ \mathbf{w}_0(
\mathbf{T}_i)\bigr\} \bigr] \bigr).
\end{eqnarray*}
%
\subsection*{Notation in step 2}
We define
$
\widehat{\mathbf{S}}_{wik}\{\bolds{\beta}_0, \widehat{\bolds
{\lambda}}(\bolds{\beta}_0, {\mathbf{w}}),\mathbf{w}(t_0)\} $ as
\begin{eqnarray*}
&& \biggl[\mathbf{Q}_{wik}\bigl\{ \widehat{\bolds{\lambda}}(\bolds{
\beta }_0, {\mathbf{w}}), \mathbf{w}(t_0)\bigr\} + {
\mathbf{Q}}_{\lambda i
k}\bigl\{ \mathbf{w}(t_0)\bigr\}^{\mathrm{T}}
\biggl\{\frac{\partial
\widehat{\bolds{\lambda}}(\bolds{\beta}_0, {\mathbf{w}})}{\partial{\mathbf{w}} } \biggr\} \biggr]
\\
&&\qquad{}\times\bigl[D_{ik} - H_{ik}\bigl\{\bolds{
\beta}_0, \widehat{\bolds{\lambda }}(\bolds{\beta}_0, {
\mathbf{w}}),\mathbf{w}(t_0)\bigr\}\bigr],
\end{eqnarray*}
and $\widehat{\mathbf{S}}_{wi}\{\bolds{\beta}_0, \widehat{\bolds
{\lambda}}(\bolds{\beta}_0, {\mathbf{w}}),\mathbf{w}(t_0)\} =
[\widehat{\mathbf{S}}_{wik}\{\bolds{\beta}_0, \widehat{\bolds
{\lambda}}(\bolds{\beta}_0, {\mathbf{w}}),\mathbf{w}(t_0)\}
^{\mathrm{T}}
, k = 1,
\ldots, M_i]^{\mathrm{T}}$.

We now define a functional from $\mathcal{D}$ to $\mathbb{R}^{d_w}$, so that
this functional evaluated at $\mathbf{w}_h$ is $\mathbf{Q}_{wik}\{
\widehat{\bolds{\lambda}}(\bolds{\beta}_0, {\mathbf{w}}),
\mathbf{w}(t_0)\} ^{\mathrm{T}}\mathbf{w}_h(t_0)$. For notational
brevity, we still use $\mathbf{Q}_{wik}\{
\widehat{\bolds{\lambda}}(\bolds{\beta}_0, {\mathbf{w}}), \mathbf{w}(t_0)\}$ to denote this functional,
that is,
\[
\mathbf{Q}_{wik}\bigl\{ \widehat{\bolds{\lambda}}(\bolds{
\beta}_0, {\mathbf{w}}), \mathbf{w}(t_0)\bigr\} (
\mathbf{w}_h) \equiv\mathbf{Q}_{wik}\bigl\{ \widehat{\bolds{
\lambda}}(\bolds{\beta}_0, {\mathbf{w}}), \mathbf{w}(t_0)
\bigr\} ^{\mathrm{T}}\mathbf{w}_h(t_0).
\]
Let\vspace*{1.5pt}
$\widehat{\mathbf{A}}_{wi}\{\bolds{\beta}_0, \widehat{\bolds
{\lambda}}(\bolds{\beta}_0, {\mathbf{w}}),\mathbf{w}(t_0)\} $ be a
$d_w \times d_w
M_i$ matrix, with the $k$th size $d_w \times d_w$ column block
$\widehat{\mathbf{A}}_{wik}\{\bolds{\beta}_0, \widehat{\bolds
{\lambda}}(\bolds{\beta}_0, {\mathbf{w}}),\mathbf{w}(t_0)\} $ being
\begin{eqnarray*}
&& \biggl[\mathbf{Q}_{wik}\bigl\{ \widehat{\bolds{\lambda}}(\bolds {
\beta}_0, {\mathbf{w}}),\mathbf{w}(t_0)\bigr\} + {
\mathbf{Q}}_{\lambda i
k}\bigl\{ \mathbf{w}(t_0)\bigr\}^{\mathrm{T}}
\biggl\{\frac{\partial
\widehat{\bolds{\lambda}}(\bolds{\beta}_0, {\mathbf{w}})}{\partial{\mathbf{w}} } \biggr\} \biggr]^{\otimes2}
\\
&&\qquad{}\times\Theta_{ik}\bigl\{\bolds{\beta}_0, \widehat{
\bolds{\lambda }}(\bolds{\beta}_0, {\mathbf{w}}),
\mathbf{w}(t_0)\bigr\}.
\end{eqnarray*}
Let\vspace*{1.5pt} $\widehat{\mathbf{V}}_{wi}\{\bolds{\beta}_0, \widehat{\bolds
{\lambda}}(\bolds{\beta}_0, {\mathbf{w}}),\mathbf{w}(t_0)\} $
be a $d_w M_i \times d_w
M_i$ matrix
with the $(p, q)$th block $\widehat{\mathbf{V}}_{wi p q}\{\bolds
{\beta}_0, \widehat{\bolds{\lambda}
}(\bolds{\beta}_0, {\mathbf{w}}),\mathbf{w}(t_0)\}$ being
\begin{eqnarray*}
&& \biggl[\mathbf{Q}_{wip}\bigl\{ \widehat{\bolds{\lambda}}(\bolds {
\beta}_0, {\mathbf{w}}), \mathbf{w}(t_0)\bigr\} + {
\mathbf{Q}}_{\lambda i
p}\bigl\{ \mathbf{w}(t_0)\bigr\}^{\mathrm{T}}
\biggl\{\frac{\partial
\widehat{\bolds{\lambda}}(\bolds{\beta}_0, {\mathbf{w}})}{\partial{\mathbf{w}} } \biggr\} \biggr]
\\
&&\qquad{}\times \biggl[\mathbf{Q}_{wiq}\bigl\{ \widehat{\bolds{\lambda}}(
\bolds {\beta}_0, {\mathbf{w}}),\mathbf{w}(t_0)\bigr\} + {
\mathbf{Q}}_{\lambda i
q}\bigl\{ \mathbf{w}(t_0)\bigr\}^{\mathrm{T}}
\biggl\{\frac{\partial
\widehat{\bolds{\lambda}}(\bolds{\beta}_0, {\mathbf{w}})}{\partial{\mathbf{w}} } \biggr\} \biggr]^{\mathrm{T}}\bolds{
\Omega}_{i pq},
\end{eqnarray*}
where $\bolds{\Omega}_{i pq}$ is the $(p,q)$th element of
the working covariance matrix $\bolds{\Omega}_{i}$.

We further define the
population level quantities ${\mathbf{S}}_{wik}\{\bolds{\beta}_0,
m_0,\mathbf{w}_0(t_0)\} $ to be
\begin{eqnarray*}
&& \bigl[\mathbf{Z}_{ik} m_0'\bigl\{
\mathbf{w}_0(t_0)^{\mathrm{T}}\mathbf{Z}_{ik}
\bigr\} - \eta\bigl\{\mathbf{w}_0(t_0)^{\mathrm{T}}
\mathbf{Z}_{ik}\bigr\} \bigr] \bigl[D_{ik} - H_{ik}
\bigl\{\bolds{\beta}_0, m_0, \mathbf{w}_0(t_0)
\bigr\}\bigr]
\end{eqnarray*}
and ${\mathbf{S}}_{wi}\{\bolds{\beta}_0, m_0,\mathbf{w}_0(t_0)\} =
[{\mathbf{S}}_{wik}\{\bolds{\beta}_0,
m_0,\mathbf{w}_0(t_0)\}^{\mathrm{T}}, k
= 1, \ldots, M_i]^{\mathrm{T}}$.
Let\break 
${\mathbf{A}}_{wi}\{\bolds{\beta}_0, m_0,\mathbf{w}_0(t_0)\} $ be a
$d_w \times d_w
M_i$ matrix, with the $k$th column block ${\mathbf{A}}_{wik}\{\bolds
{\beta}_0, m_0,\mathbf{w}
_0(t_0)\}$ being a $d_w \times d_w$ matrix
\begin{eqnarray*}
&& \bigl[\mathbf{Z}_{ik} m_0'\bigl\{
\mathbf{w}_0(t_0)^{\mathrm
{T}}\mathbf{Z}_{ik}
\bigr\} - \eta\bigl\{\mathbf{w}_0(t_0)^{\mathrm{T}}
\mathbf{Z}_{ik}\bigr\} \bigr]^{\otimes2} \Theta_{ik}\bigl
\{\bolds{\beta}_0, m_0, \mathbf{w}_0(t_0)
\bigr\}.
\end{eqnarray*}
Let ${\mathbf{V}}_{wi}\{\bolds{\beta}_0, m_0,\mathbf{w}_0(t_0)\} $
be a $d_w M_i \times d_w
M_i$ matrix
with the $(p, q)$th block
${\mathbf{V}}_{wi pq}\{\bolds{\beta}_0, m_0, \mathbf{w}_0(t_0)\}$ being
\begin{eqnarray*}
&&\bigl[\mathbf{Z}_{ip} m'_0\bigl\{
\mathbf{w}_0(t_0)^{\mathrm{T}}\mathbf{Z}_{ip}
\bigr\} - \eta\bigl\{\mathbf{w}_0(t_0)^{\mathrm{T}}
\mathbf{Z}_{ip}\bigr\} \bigr]
\\
&&\qquad{}\times\bigl[\mathbf{Z}_{iq} m'_0\bigl\{
\mathbf{w}_0(t_0)^{\mathrm
{T}}\mathbf{Z}_{iq}
\bigr\} - \eta\bigl\{\mathbf{w} _0(t_0)^{\mathrm{T}}
\mathbf{Z}_{iq}\bigr\} \bigr]^{\mathrm{T}}\bolds{\Omega
}_{i pq}.
\end{eqnarray*}
Let
${\mathbf{V}^*}_{wi}\{\bolds{\beta}_0, m_0,\mathbf{w}_0(t_0)\} $ be
a $d_w M_i \times d_w
M_i$ matrix. The $(p, q)$th block is obtained by
replacing $\bolds{\Omega}_{i pq}$ in ${\mathbf{V}}_{wi pq}\{\bolds
{\beta}_0, m_0, \mathbf{w}_0(t_0)\}$
with
\[
\bigl[E(D_{ip}D_{iq}) - H_{ip}\bigl\{\bolds{
\beta}_0, m_0,\mathbf{w}_0(t_0)
\bigr\} H_{iq}\bigl\{\bolds {\beta}_0, m_0,
\mathbf{w}_0(t_0)\bigr\}\bigr].
\]

Here $\eta$ is an operator that maps functions in $C^1([0, \tau])$ to
functionals from $\mathcal{D}$ to~$\mathbb{R}^{d_w}$.
Specifically, $\eta$ minimizes
\begin{eqnarray*}
&&\sup_{\mathbf{w}_{h}\in\mathcal{D}} \bigl\|E \bigl( \bigl[\widetilde{
\mathbf{Q}}_{\mathbf{w}i}\bigl\{m_0, \mathbf{w} _h(
\mathbf{T}_i)\bigr\}- \eta\bigl\{\mathbf{U}_i(
\mathbf{T}_i)\bigr\}(\mathbf{w}_h) \bigr]^{\mathrm
{T}}
\Theta_i \bigl\{\bolds{\beta}_0, m_0,
\mathbf{w}_0(\mathbf{T}_{i})\bigr\}
\\
&&\hspace*{16pt}\qquad{}\times\bolds{\Omega}_i^{-1} \Theta_i \bigl
\{\bolds{\beta}_0, m_0,\mathbf{w}_0(
\mathbf{T}_{i})\bigr\} \bigl[\widetilde{\mathbf{Q}}_{\mathbf{w}i}
\bigl\{m_0, \mathbf{w} _h(\mathbf{T}_i)\bigr
\}- \eta\bigl\{\mathbf{U}_i(\mathbf{T}_i)\bigr\}(
\mathbf{w}_h) \bigr] \bigr) \bigr\|_2,
\end{eqnarray*}
where
\begin{eqnarray*}
\widetilde{\mathbf{Q}}_{wi}\bigl\{m_0,
\mathbf{w}_h(\mathbf{T}_i)\bigr\} &=&
\bigl[m'_0\bigl\{\mathbf{w}_0(T_{i1})^{\mathrm{T}}
\mathbf{Z}_{i1}\bigr\}\mathbf{w}_h(T_{i1})^{\mathrm{T}}
\mathbf{Z}_{i1}, \ldots,
\\
&&\hspace*{6pt}m'_0\bigl\{\mathbf{w}_0(T_{iM_i})^{\mathrm{T}}
\mathbf{Z}_{iM_i}\bigr\} \mathbf{w}_h(T_{iM_i})^{\mathrm{T}}
\mathbf{Z}_{i
M_i}\bigr]^{\mathrm{T}}
\end{eqnarray*}
and
$
\eta\{\mathbf{U}_i(\mathbf{T}_i)\}(\mathbf{w}_h)$ $=$ $
[\eta\{\mathbf{w}(T_{ik})^{\mathrm{T}}\mathbf{Z}_{ik}\}(\mathbf{w}_h), k = 1, \dots,
M_i]^{\mathrm{T}}$ are $M_i$ vectors. We can also write
\[
\eta\bigl(\mathbf{w}_0(T_{ik})^{\mathrm{T}}
\mathbf{Z}_{ik}\bigr) = E\bigl[\mathbf{Z}_{ik}
m'_0\bigl\{\mathbf{w}_0(T_{ik})^{\mathrm{T}}
\mathbf{Z}_{ik}\bigr\}| \mathbf{w}_0(T_{ik})^{\mathrm{T}}
\mathbf{Z}_{ik}\bigr].
\]
Further,\vspace*{1.5pt} we define $\widetilde{\mathbf{Q}}_{wi}\{\widehat{\bolds
{\lambda}}\{\widehat{\bolds{\beta}}, \widehat{\mathbf{w}
}(\widehat{\bolds{\beta}})\}, \cdot\}$ as a ${M_i}\times d_w$
matrix, with row $j$ as
$\mathbf{B}_r'\{\widehat{\mathbf{w}}(\widehat{\bolds{\beta}},
T_{ik})^{\mathrm{T}}\mathbf{Z}_{ik}\}\widehat{\bolds{\lambda}}\{
\widehat{\bolds{\beta}}, \widehat{\mathbf{w}}(\widehat{\bolds
{\beta}})\}\mathbf{Z}_{ik}^{\mathrm{T}}$.
In the estimation, we use the asymptotic form in Lemma~4
in the supplementary article in place of $\partial
\widehat{\bolds{\lambda}}(\bolds{\beta}_0, {\mathbf{w}})/\partial
{\mathbf{w}} $ for computation.

\subsection*{Notation in step 3}
We define
\begin{eqnarray*}
&&\widehat{\mathbf{S}}_{\beta i k}\bigl[\bolds{\beta}, \widehat{\bolds {
\lambda}}\bigl\{\bolds{\beta}, \widehat{\mathbf{w}}(\bolds{\beta},
T_{ik})\bigr\}, {\widehat{\mathbf{w}}}(\bolds{\beta})\bigr]
\\
&&\qquad= \biggl(\mathbf{Q}_{\beta i k} + \biggl[\frac{\partial\widehat
{\bolds{\lambda}}\{\bolds{\beta},
{\widehat{\mathbf{w}}}(\bolds{\beta})\}
}{\partial\bolds{\beta}^{\mathrm{T}}} +
\frac{\partial\widehat
{\bolds{\lambda}}\{\bolds{\beta}, {\widehat{\mathbf{w}
}}(\bolds{\beta})\}
}{ {\widehat{\mathbf{w}}}(\bolds{\beta})}\frac{\partial\widehat
{\mathbf{w}}(\bolds{\beta})}{\partial
\bolds{\beta}^{\mathrm{T}}} \biggr]^{\mathrm{T}}\mathbf{Q}_{\lambda
i k}
\bigl\{ \widehat{\mathbf{w}}(\bolds{\beta}, T_{ik})\bigr\}
\\
&&\hspace*{80pt}\qquad\quad{}+ \biggl\{\frac{\partial\widehat{\mathbf{w}}(\bolds{\beta},
T_{ik})}{\partial
\bolds{\beta}^{\mathrm{T}}} \biggr\}^{\mathrm{T}}\mathbf{Q}_{w i k}
\bigl[ \widehat{\bolds{\lambda}}\bigl\{\bolds{\beta}, {\widehat{\mathbf{w}}}(
\bolds{\beta})\bigr\}, \widehat{\mathbf{w}}(\bolds{\beta}, T_{ik})\bigr]
\biggr)
\\
&&\qquad\quad{}\times \bigl(D_{ik} - H_{ik} \bigl[\bolds{\beta},
\widehat{\bolds {\lambda}}\bigl\{\bolds{\beta}, {\widehat{\mathbf{w}}}(\bolds{
\beta})\bigr\}, \widehat{\mathbf{w}}(\bolds{\beta}, T_{ik})\bigr]
\bigr),
\end{eqnarray*}
and $\widehat{\mathbf{S}}_{\beta i}[\bolds{\beta}, \widehat{\bolds
{\lambda}}\{\bolds{\beta}, {\widehat{\mathbf{w}}}(\bolds{\beta
})\},
\widehat{\mathbf{w}}(\bolds{\beta}, \mathbf{T}_{i})] = (\widehat
{\mathbf{S}}_{\beta i k}[\bolds{\beta}, \widehat{\bolds{\lambda
}}\{
\bolds{\beta}, \widehat{\mathbf{w}}(\bolds{\beta}, T_{ik})\}
^{\mathrm{T}},
{\widehat{\mathbf{w}}}(\bolds{\beta})],\break  k = 1, \ldots,
M_i)^{\mathrm{T}}$.
Let $\widehat{\mathbf{A}}_{\beta i}[\bolds{\beta}, \widehat{\bolds
{\lambda}}\{\bolds{\beta}, \widehat{\mathbf{w}}(\bolds{\beta},
\mathbf{T}
_{i})\},
{\widehat{\mathbf{w}}}(\bolds{\beta})]$ be a $d_\beta\times
d_\beta M_i$ matrix with the $k$th
size $d_\beta\times d_\beta$ column block $\widehat{\mathbf{A}}_{\beta i k}[\bolds{\beta}
, \widehat{\bolds{\lambda}}\{\bolds{\beta}, {\widehat{\mathbf
{w}}}(\bolds{\beta})\},
\widehat{\mathbf{w}}(\bolds{\beta}, T_{ik})]$ being
\begin{eqnarray*}
&& \biggl(\mathbf{Q}_{\beta i k}+ \biggl[\frac{\partial\widehat
{\bolds{\lambda}}\{\bolds{\beta},
{\widehat{\mathbf{w}}}(\bolds{\beta})\}
}{\partial\bolds{\beta}^{\mathrm{T}}} +
\frac{\partial\widehat
{\bolds{\lambda}}\{\bolds{\beta}, {\widehat{\mathbf{w}
}}(\bolds{\beta})\}
}{ {\widehat{\mathbf{w}}}(\bolds{\beta})}\frac{\partial\widehat
{\mathbf{w}}(\bolds{\beta})}{\partial
\bolds{\beta}^{\mathrm{T}}} \biggr]^{\mathrm{T}}\mathbf{Q}_{\lambda
i k}
\bigl\{ \widehat{\mathbf{w}}(\bolds{\beta}, T_{ik})\bigr\}
\\
&&\hspace*{58pt}\qquad {}+ \biggl\{\frac{\partial\widehat{\mathbf{w}}(\bolds{\beta},
T_{ik})}{\partial
\bolds{\beta}^{\mathrm{T}}} \biggr\}^{\mathrm{T}}\mathbf{Q}_{w i k}
\bigl[ \widehat{\bolds{\lambda}}\bigl\{\bolds{\beta}, {\widehat{\mathbf{w}}}(
\bolds{\beta})\bigr\}, \widehat{\mathbf{w}}(\bolds{\beta}, T_{ik})\bigr]
\biggr)^{\otimes2} \\
&&\qquad{}\times\bolds{\Theta}_i \bigl[\bolds{\beta},
\widehat{\bolds{\lambda}}\bigl\{\bolds{\beta}, {\widehat{\mathbf{w}}}(\bolds{
\beta})\bigr\},
 \widehat{\mathbf{w}}(\bolds{\beta}, T_{ik})\bigr].
\end{eqnarray*}
Let\vspace*{1.5pt} $\widehat{\mathbf{V}}_{\beta
i}[\bolds{\beta}, \widehat{\bolds{\lambda}}\{\bolds{\beta},
{\widehat{\mathbf{w}}}(\bolds{\beta})\}, \widehat{\mathbf{w}}(\bolds{\beta},
\mathbf{T}_{i})]^{-1} $ be a $d_\beta M_i \times d_\beta M_i$ matrix
with the
$(p, q)$th
block $\widehat{\mathbf{V}}_{\beta i p }[\bolds{\beta}, \widehat
{\bolds{\lambda}}\{\bolds{\beta}, {\widehat{\mathbf{w}
}}(\bolds{\beta})\},
\widehat{\mathbf{w}}(\bolds{\beta}, T_{ip})] $ being
\begin{eqnarray*}
&& \biggl(\mathbf{Q}_{\beta i p} + \biggl[\frac{\partial\widehat
{\bolds{\lambda}}\{\bolds{\beta},
{\widehat{\mathbf{w}}}(\bolds{\beta})\}
}{\partial\bolds{\beta}^{\mathrm{T}}} +
\frac{\partial\widehat
{\bolds{\lambda}}\{\bolds{\beta}, {\widehat{\mathbf{w}
}}(\bolds{\beta})\}
}{ {\widehat{\mathbf{w}}}(\bolds{\beta})}\frac{\partial\widehat
{\mathbf{w}}(\bolds{\beta})}{\partial
\bolds{\beta}^{\mathrm{T}}} \biggr]^{\mathrm{T}}\mathbf{Q}_{\lambda
i p}
\bigl\{ \widehat{\mathbf{w}}(\bolds{\beta}, T_{ip})\bigr\}
\\
&&\hspace*{82pt}{}+ \biggl\{\frac{\partial\widehat{\mathbf{w}}(\bolds{\beta},
T_{ip})}{\partial
\bolds{\beta}^{\mathrm{T}}} \biggr\}^{\mathrm{T}}\mathbf{Q}_{w i p}
\bigl[ \widehat{\bolds{\lambda}}\bigl\{\bolds{\beta}, {\widehat{\mathbf{w}}}(
\bolds{\beta})\bigr\}, \widehat{\mathbf{w}}(\bolds{\beta}, T_{ip})\bigr]
\biggr)
\\
&&\qquad{}\times\biggl(\mathbf{Q}_{\beta iq} + \biggl[\frac
{\partial\widehat{\bolds{\lambda}}\{\bolds{\beta}, {\widehat
{\mathbf{w}}}(\bolds{\beta})\}
}{\partial\bolds{\beta}^{\mathrm{T}}}+ \frac{\partial\widehat{\bolds{\lambda}}\{\bolds{\beta},
{\widehat{\mathbf{w}}}(\bolds{\beta})\}
}{ {\widehat{\mathbf{w}}}(\bolds{\beta})}\frac{\partial\widehat
{\mathbf{w}}(\bolds{\beta})}{\partial
\bolds{\beta}^{\mathrm{T}}} \biggr]^{\mathrm{T}}
\mathbf{Q}_{\lambda
iq} \bigl\{ \widehat{\mathbf{w}}(\bolds{\beta},
T_{iq})\bigr\}
\\[1pt]
&&\hspace*{92pt}\qquad{}+ \biggl\{\frac{\partial\widehat{\mathbf{w}}(\bolds{\beta},
T_{iq})}{\partial
\bolds{\beta}^{\mathrm{T}}} \biggr\}^{\mathrm{T}}\mathbf{Q}_{w iq} \bigl[\widehat{\bolds{\lambda}}\bigl\{\bolds
{\beta}, {\widehat{\mathbf{w}}}(\bolds{\beta})\bigr\}, \widehat{\mathbf{w}}(
\bolds{\beta}, T_{iq})\bigr] \biggr)^{\mathrm
{T}}\\[1pt]
&&\qquad{}\times\bolds{
\Omega}_{ip q}.
\end{eqnarray*}
Additionally,\vspace*{1pt} let $\delta_u \in
C^q([0, 1])$, and we define $\bolds{\delta}\{\mathbf{w}(T_{ik})^{\mathrm{T}}
\mathbf{Z}_{ik}\} =  \break [\delta_u
\{\mathbf{w}(T_{ik})^{\mathrm{T}}\mathbf{Z}_{ik}\}$, $u = 1, \ldots
, d_\beta]\in\mathrm
{R}^{d_\beta}$
which minimizes
\begin{eqnarray*}
&&{\mathbf1}_{d_\beta}^{\mathrm{T}}E \bigl( \bigl[\widetilde{
\mathbf{Q}}_{\beta i}- \bolds{\delta}\bigl\{\mathbf{U}_i(
\mathbf{T}_i)\bigr\} \bigr]^{\mathrm{T}}\Theta _i \bigl
\{\bolds{\beta}_0, m_0,\mathbf{w}_0(\mathbf{T}
_{i})\bigr\} \bolds{\Omega}_i^{-1}
\Theta_i \bigl\{\bolds{\beta}_0, m_0,
\mathbf{w}_0(\mathbf{T}_{i})\bigr\}
\\[1pt]
&& \hspace*{218pt}{}\times\bigl[\widetilde{\mathbf{Q}}_{\beta i}- \bolds{\delta}\bigl\{
\mathbf{U}_i(\mathbf{T}_i)\bigr\} \bigr] \bigr){
\mathbf1}_{d_\beta},
\end{eqnarray*}
where\vspace*{1pt}
$
\widetilde{\mathbf{Q}}_{\beta i} =(\mathbf{X}_{i1}, \ldots, \mathbf{X}_{i M_i})^{\mathrm{T}}
$
is a $M_i \times d_\beta$ matrix, and
$
\bolds{\delta}\{\mathbf{U}_i(\mathbf{T}_i)\}=\break 
[\bolds{\delta}\{\mathbf{w}(T_{ik})^{\mathrm{T}} \mathbf{Z}_{ik}\},
k = 1, \dots, M_i]^{\mathrm{T}}$
is a $M_i \times d_\beta$
matrix. We can\vspace*{1pt} also write
$\bolds{\delta}\{\mathbf{w}_0(T_{ik})^{\mathrm{T}}\mathbf{Z}_{ik}\}$
as $E\{\mathbf{X}
|\mathbf{w}_0(T_{ik})^{\mathrm{T}}\mathbf{Z}_{ik} \}$.
Further, we define
\begin{eqnarray*}
{\mathbf{B}}(t_0) &=& E \bigl( {\mathbf{A}}_{wi}\bigl\{
\bolds{\beta}_0, m_0,\mathbf{w}_0(t_0)
\bigr\} {\mathbf{V}}_{wi}\bigl\{\bolds{\beta}_0,
m_0,\mathbf{w}_0(t_0)\bigr\}^{-1}\\[1pt]
&&\quad{}\times
\bigl[\mathbf{Q}_{wi}\bigl\{ m_0,\mathbf{w}_0(t_0)
\bigr\}- \eta\bigl\{\mathbf{U}_i(t_0)\bigr\} \bigr]
\nonumber
\\[1pt]
&&\hspace*{19pt}\quad{}\times\bolds{
\Theta}^*_i \bigl\{\bolds{\beta }_0, m_0,
\mathbf{w}_0(t_0)\bigr\} \mathbf{Q}^*_{w i}\bigl
\{ m_0,\mathbf{w}_0(t_0)\bigr\} \bigr),
\end{eqnarray*}
where
$\bolds{\Theta}^*_i \{\bolds{\beta}_0, m_0,\mathbf{w}(t_0)\} $
is a $d_w M_i \times d_w M_i$ diagonal matrix with the $k$th diagonal
block being a $d_w \times d_w$ diagonal with the element
$\Theta_{ik}\{\bolds{\beta}_0, m_0,\mathbf{w}(t_0)\}$,
and
$\mathbf{Q}_{wi}\{m_0, \mathbf{w}(t_0)\}$
is a $d_w {M_i} \times d_w {M_i}$ diagonal\vspace*{1pt} matrix with the $k$th
diagonal block being\vspace*{1pt} $\operatorname{diag}
[\mathbf{Z}_{ik} m'_0\{\mathbf{w}(t_0)^{\mathrm{T}}\mathbf{Z}_{ik}\}]$.
Moreover
$\mathbf{Q}^*_{wi}\{m_0, \mathbf{w}(\mathbf{t}_0)\}
$
is a $d_w M_i \times d_w $ matrix with the $k$th row block being a $d_w
\times d_w$ matrix with $d_w$\vspace*{1pt} replications of $\mathbf
{Z}_{ik}^{\mathrm{T}}
m'_0\{\mathbf{w}(t_0)^{\mathrm{T}}\mathbf{Z}_{ik}\}$, and\vspace*{1pt} $\eta\{
\mathbf{U}_i(t_0)\} = [\eta
\{\mathbf{w}(t_0)\mathbf{Z}_{i1}\}^{\mathrm{T}}, \ldots,\eta
\{\mathbf{w}(t_0)\mathbf{Z}_{iM_i}\}^{\mathrm{T}}]^{\mathrm{T}}$.\vspace*{1pt}
Also let $\mathbf{B}(\mathbf{T}_i)$ be the $d_w
M_i\times d_w M_i$ block diagonal\vspace*{1pt} matrix with the $k$th block as
$\mathbf{B}(T_{ik})$ and $\mathbf{f}_\mathbf{T}(\mathbf{T}_i)$ be the $d_w
M_i\times d_w M_i$ block diagonal matrix with\vspace*{1pt} the $k$th block as
$f_T(T_{ik})$.

Let\vspace*{1pt} $\gamma_u \in
C^q([0, 1])$, and we define $\bolds{\gamma}\{\mathbf
{w}(T_{ik})^{\mathrm{T}}
\mathbf{Z}_{ik}\} = [\gamma_u
\{\mathbf{w}(T_{ik})^{\mathrm{T}}\mathbf{Z}_{ik}\}$, $u = 1, \ldots
, d_\beta]\in
\mathrm{R}^{d_\beta}$, which minimize
\begin{eqnarray*}
&&{\mathbf1}_{\beta}^{\mathrm{T}}E \biggl[ \biggl\{\widetilde{\mathbf
{Q}}_{\mathbf{w}i} \biggl(m_0, \mathbf{B} (\mathbf{T}_{i})^{-1}E
\biggl[{\mathbf{A}}_{wj}\bigl\{\bolds{\beta}_0,
m_0,\mathbf{w}_0(\mathbf{T}_{i})\bigr\}{
\mathbf{V}}_{wj}\bigl\{\bolds{\beta }_0, m_0,\mathbf{w}_0(\mathbf{T}_{i})\bigr\}^{-1}
\\[1pt]
&&\hspace*{153pt}\qquad\qquad{}\times
\frac{\partial{\mathbf
{S}}_{wj}\{\bolds{\beta}_0,
m_0,\mathbf{w}_0(\mathbf{T}_{i})\} }{\partial\bolds{\beta
}^{\mathrm{T}}} \Big| \mathbf{O}_i \biggr]  \biggr)\\[1pt]
 &&\hspace*{231pt}\qquad\qquad{}- \bolds{\gamma}
\bigl\{\mathbf{U}_i(\mathbf{T}_i)\bigr\} \biggr\}^{\mathrm
{T}}
\\[1pt]
&&\hspace*{24pt}{}\times\Theta_i \bigl\{\bolds{\beta}_0,
m_0,\mathbf{w}_0(\mathbf{T}_{i})\bigr\} \bolds{
\Omega}_i^{-1} \Theta_i \bigl\{\bolds{
\beta}_0, m_0,\mathbf{w}_0(
\mathbf{T}_{i})\bigr\}\\
&&\hspace*{24pt}{}\times \biggl\{\widetilde{\mathbf{Q}}_{\mathbf{w}i}
\biggl(m_0, \mathbf{B}(\mathbf{T}_{i})^{-1}E
\biggl[{\mathbf{A}}_{wj}\bigl\{ \bolds{\beta}_0,
m_0,\mathbf{w}_0(\mathbf{T}_{i})\bigr\}{
\mathbf{V}}_{wj}\bigl\{\bolds{\beta }_0, m_0,
\mathbf{w}_0(\mathbf{T}_{i})\bigr\}^{-1}
\\[-2pt]
&&\hspace*{165pt}\qquad\qquad{}\times
\frac
{\partial{\mathbf{S}}_{wj}\{\bolds{\beta}_0,
m_0,\mathbf{w}_0(\mathbf{T}_{i})\} }{\partial\bolds{\beta
}^{\mathrm{T}}}  \Big| \mathbf{O}_i \biggr]\biggr)
\\[-2pt]
&& \hspace*{269pt}{}- \bolds{\gamma}\bigl\{\mathbf{U}_i(\mathbf{T}_i)
\bigr\} \biggr\} \biggr]\mathbf{1}_{\beta},
\end{eqnarray*}
where $
\bolds{\gamma}\{\mathbf{U}_i(\mathbf{T}_i)\}$ $=$ $
[\bolds{\gamma}\{\mathbf{w}(T_{ik})^{\mathrm{T}}\mathbf{Z}_{ik}\},
k = 1, \dots, M_i]^{\mathrm{T}}$ is a $M_i
\times d_\beta$,\vspace*{-2pt} and
\begin{eqnarray*}
&&\widetilde{\mathbf{Q}}_{wi} \biggl(m_0, \mathbf{B}(
\mathbf {T}_{i})^{-1}E \biggl[{\mathbf{A}}_{wj}
\bigl\{\bolds{\beta}_0, m_0,\mathbf{w}_0(
\mathbf{T}_{i})\bigr\}{\mathbf{V}}_{wj}\bigl\{\bolds{\beta
}_0, m_0,\mathbf{w}_0(\mathbf{T}_{i})
\bigr\}^{-1}
\\[-2pt]
&&\hspace*{166pt}{}\times\frac{\partial{\mathbf{S}}_{wj}\{\bolds{\beta}_0,
m_0,\mathbf{w}_0(\mathbf{T}_{i})\} }{\partial\bolds{\beta
}^{\mathrm{T}}} \Big| \mathbf{O}_i \biggr]  \biggr)
\end{eqnarray*}
is a $M_i\times\bolds{\beta}$ matrix with $k$th row\vspace*{-2pt} as
\begin{eqnarray*}
&& \biggl({\mathbf{B}}(T_{ik})^{-1} E \biggl[{
\mathbf{A}}_{wj}\bigl\{\bolds {\beta}_0, m_0,
\mathbf{w}_0(T_{ik})\bigr\}{\mathbf{V}}_{wj}\bigl
\{\bolds{\beta}_0, m_0,\mathbf{w}_0(T_{ik})
\bigr\}^{-1}
\\[-2pt]
&&\hspace*{115pt}\qquad{}\times\frac{\partial{\mathbf{S}}_{wj}\{\bolds{\beta}_0,
m_0,\mathbf{w}_0(T_{ik})\} }{\partial\bolds{\beta}^{\mathrm
{T}}}\Big\vert \mathbf{O}_i \biggr]
\biggr)^{\mathrm{T}}\\
&&\qquad{}\times\mathbf{Z}_{ik} m'_0\bigl
\{\mathbf{w}_0(T_{ik})^{\mathrm{T}}\mathbf{Z}_{ik}
\bigr\}.
\end{eqnarray*}
We can also write
\begin{eqnarray*}
&& \bolds{\gamma}\bigl(\mathbf{w}_0(T_{ik})^{\mathrm{T}}
\mathbf{Z}_{ik}\bigr) \\
&&\qquad= E \biggl\{ \biggl({\mathbf{B}}(T_{ik})^{-1}
E \biggl[{\mathbf {A}}_{wj}\bigl\{\bolds{\beta}_0,
m_0,\mathbf{w}_0(T_{ik})\bigr\}{
\mathbf{V}}_{wj}\bigl\{\bolds{\beta}_0, m_0,\mathbf{w}_0(T_{ik})\bigr\}^{-1}
\\[-2pt]
&&\hspace*{184pt}{}\times
\frac{\partial{\mathbf{S}}_{wj}\{
\bolds{\beta}_0,
m_0,\mathbf{w}_0(T_{ik})\} }{\partial\bolds{\beta}^{\mathrm
{T}}}\Big\vert \mathbf{O}_i \biggr] \biggr)^{\mathrm{T}}
\mathbf{Z}_{ik}
\\[-2pt]
&&\hspace*{197pt}{}\times m'_0\bigl\{\mathbf{w}_0(T_{ik})^{\mathrm{T}}
\mathbf{Z}_{ik}\bigr\} | \mathbf{w}_0(T_{ik})^{\mathrm{T}}
\mathbf{Z} _{ik} \biggr\}.
\end{eqnarray*}
We also define the population forms ${\mathbf{S}}_{\beta i k}\{\bolds
{\beta}_0,
m_0, {\mathbf{w}}_0(T_{ik})\}$ as
\begin{eqnarray*}
&& \biggl\{\mathbf{Q}_{\beta i k} - \bolds{\delta}\bigl\{\mathbf
{w}_0(T_{ik})^{\mathrm{T}}\mathbf{Z}_{ik}\bigr\}
\\[-2pt]
&&\hspace*{6pt}{}- \biggl({\mathbf{B}}(T_{ik})^{-1} E \biggl[{
\mathbf{A}}_{wj}\bigl\{\bolds {\beta}_0, m_0,
\mathbf{w}_0(T_{ik})\bigr\}
\\[-2pt]
&&
\hspace*{34pt}\qquad\qquad{}\times{\mathbf{V}}_{wj}\bigl\{\bolds{\beta}_0,
m_0,\mathbf{w}_0(T_{ik})\bigr\} ^{-1}\frac{\partial{\mathbf{S}
}_{wj}\{\bolds{\beta}_0,
m_0,\mathbf{w}_0(T_{ik})\} }{\partial\bolds{\beta}^{\mathrm
{T}}}\Big\vert \mathbf{O}_i \biggr] \biggr)^{\mathrm{T}}
\mathbf{Z}_{ik}
\\[-2pt]
&&\hspace*{130pt}\qquad\qquad{}\times m'_0\bigl\{\mathbf{w}_0(T_{ik})^{\mathrm{T}}
\mathbf{Z}_{ik}\bigr\}+ \bolds {\gamma}\bigl\{\mathbf{w}_0(T_{ik})^{\mathrm{T}}
\mathbf{Z} _{ik}\bigr\} \biggr\} \\
&&\qquad{}\times\bigl[D_{ik} -
H_{ik} \bigl\{\bolds{\beta}_0, m_0, {
\mathbf{w}}_0(T_{ik})\bigr\} \bigr]
\end{eqnarray*}
and ${\mathbf{S}}_{\beta i}\{\bolds{\beta}_0, m_0, {\mathbf
{w}}_0(\mathbf{T}_{i})\} = [{\mathbf{S}}_{\beta i
k}\{\bolds{\beta}_0,
m_0, {\mathbf{w}}_0(T_{ik})\}^{\mathrm{T}}, k = 1, \ldots,
M_i]^{\mathrm{T}}$.
Let\break ${\mathbf{A}}_{\beta i}\{\bolds{\beta}_0, m_0, {\mathbf
{w}}_0(\mathbf{T}_{i})\}$ be
a $d_\beta\times d_\beta M_i $ be the matrix with the $k$th block
${\mathbf{A}}_{\beta ik}\{\bolds{\beta}_0, m_0, {\mathbf
{w}}_0(T_{ik})\} $ being a
$d_\beta\times d_\beta$ matrix
\begin{eqnarray*}
&& \biggl\{\mathbf{Q}_{\beta i k} - \bolds{\delta}\bigl\{\mathbf
{w}_0(T_{ik})^{\mathrm{T}}\mathbf{Z}_{ik}\bigr\}
\\
&&\hspace*{6pt}{}- \biggl({\mathbf{B}}(T_{ik})^{-1} E \biggl[{
\mathbf{A}}_{wj}\bigl\{\bolds {\beta}_0, m_0,
\mathbf{w}_0(T_{ik})\bigr\}
\\
&&\hspace*{56pt}\qquad{}\times{\mathbf{V}}_{wj}\bigl\{\bolds{\beta}_0,
m_0,\mathbf{w}_0(T_{ik})\bigr\} ^{-1}
\\
&&\hspace*{56pt}\qquad{}\times\frac{\partial{\mathbf{S}
}_{wj}\{\bolds{\beta}_0,
m_0,\mathbf{w}_0(T_{ik})\} }{\partial\bolds{\beta}^{\mathrm
{T}}}\Big\vert \mathbf{O}_i \biggr] \biggr)^{\mathrm{T}}
\mathbf{Z}_{ik}
\\
&&\qquad\qquad\hspace*{21pt}{}\times m'_0\bigl\{\mathbf{w}_0(T_{ik})^{\mathrm{T}}
\mathbf{Z}_{ik}\bigr\} + \bolds{\gamma}\bigl\{\mathbf{w}
_0(T_{ik})^{\mathrm{T}}\mathbf{Z}_{ik}\bigr\}
\biggr\} ^{\otimes2}\\
&&\qquad{}\times \Theta _{ik}\bigl[\bolds{\beta}_0,
m_0,\mathbf{w}_0(T_{ik})\bigr].
\end{eqnarray*}
Let
${\mathbf{V}}_{\beta i}\{\bolds{\beta}_0, m_0,{\mathbf
{w}}_0(\mathbf{T}_{i})\} $ be a $d_\beta M_i
\times d_\beta M_i$ with the $(p, q)$th block ${\mathbf{V}}_{\beta i p
q}\{\bolds{\beta}_0,\break  m_0, {\mathbf{w}}_0(T_{ip})\}$ being
\begin{eqnarray*}
&&\!\!\! \biggl\{\mathbf{Q}_{\beta i p} - \bolds{\delta}\bigl\{\mathbf
{w}_0(T_{ip})^{\mathrm{T}}\mathbf{Z}_{ip}\bigr\}
\\
&&\!\!\!\hspace*{6pt}{}- \biggl(\mathbf{B}(T_{ip})^{-1} E \biggl[{
\mathbf{A}}_{wj}\bigl\{\bolds {\beta}_0, m_0,
\mathbf{w}_0(T_{ip})\bigr\}
\\
&&\!\!\!\hspace*{80pt}{}\times{\mathbf{V}}_{wj}\bigl\{\bolds{\beta}_0,
m_0,\mathbf{w}_0(T_{ip})\bigr\} ^{-1}
\\
&&\!\!\!\hspace*{80pt}{}\times\frac{\partial{\mathbf{S}
}_{wj}\{\bolds{\beta}_0,
m_0,\mathbf{w}_0(T_{ip})\} }{\partial\bolds{\beta}^{\mathrm
{T}}}\Big\vert \mathbf{O}_i \biggr] \biggr)^{\mathrm{T}}
\mathbf{Z}_{ip}
\\
&&\!\!\!\hspace*{153pt}{}\times m'_0\bigl\{\mathbf{w}_0(T_{ip})^{\mathrm{T}}
\mathbf{Z}_{ip}\bigr\}\biggr\}\\
&&\!\!\!\qquad{}+ \bolds{\gamma}\bigl\{
\mathbf{w}_0(T_{ip})^{\mathrm{T}}\mathbf{Z}_{ip}
\bigr\}  \biggl\{\mathbf{Q}_{\beta i q} - \bolds{\delta}\bigl\{
\mathbf{w}_0(T_{iq})^{\mathrm{T}}\mathbf{Z}_{iq}
\bigr\}
\\
&&\!\!\! \hspace*{112pt}{}- \biggl(\mathbf{B}(T_{iq})^{-1} E \biggl[{
\mathbf{A}}_{wj}\bigl\{\bolds {\beta}_0, m_0,
\mathbf{w}_0(T_{iq})\bigr\}\\
&&\!\!\!\hspace*{185pt}{}\times{\mathbf{V}}_{wj}\bigl
\{\bolds{\beta}_0, m_0,\mathbf{w}_0(T_{iq})
\bigr\}^{-1}
\\
&&\!\!\!\hspace*{185pt}{}\times\frac{\partial{\mathbf{S}}_{wj}\{\bolds{\beta}_0,
m_0,\mathbf{w}_0(T_{iq})\} }{\partial\bolds{\beta}^{\mathrm
{T}}}\Big\vert \mathbf{O}_i \biggr]
\biggr)^{\mathrm{T}}\\
&&\!\!\!\hspace*{130pt}{}\times\mathbf{Z}_{iq} m'_0\bigl
\{\mathbf{w}_0(T_{iq})^{\mathrm{T}}\mathbf{Z}_{iq}
\bigr\}
+ \bolds{\gamma}\bigl\{\mathbf{w}_0(T_{iq})^{\mathrm{T}}
\mathbf {Z}_{iq}\bigr\} \biggr\} ^{\mathrm{T}}
\bolds{
\Omega}_{i p q}.
\end{eqnarray*}
Let ${\mathbf{V}}^*_{\beta i}\{\bolds{\beta}_0, m_0,{\mathbf
{w}}_0(\mathbf{T}_{i})\} $ be a $d_\beta
M_i \times d_\beta
M_i$ matrix. The $(p, q)$th block is obtained by
replacing $\bolds{\Omega}_{i p q}$ in ${\mathbf{V}}_{\beta i}\{
\bolds{\beta}_0,
m_0,{\mathbf{w}}_0(\mathbf{T}_{i})\} $ with
\[
\bigl[E(D_{ip}D_{iq}) - H_{ip}\bigl\{\bolds{
\beta}_0, m_0,\mathbf{w}(T_{ip})\bigr\}
H_{iq}\bigl\{\bolds {\beta}_0, m_0,
\mathbf{w}(T_{iq})\bigr\}\bigr].
\]
\end{appendix}

\section*{Acknowledgments}

The authors thank the Editor, Associate Editor and
three anonymous referees for their comprehensive reviews, which greatly
improved the paper.

\begin{supplement}[id=suppA]
\stitle{Supplement to ``Fused kernel-spline smoothing for repeatedly measured
outcomes in a generalized partially linear model with functional single
index''}
\slink[doi]{10.1214/15-AOS1330SUPP} 
\sdatatype{.pdf}
\sfilename{aos1330\_supp.pdf}
\sdescription{We provide\vadjust{\goodbreak} the comprehensive proofs of Theorems~\ref{th1},
\ref{th2} and \ref{th3} and additional lemmas which support the results.}
\end{supplement}


%


%

%


\printaddresses
\end{document}